\documentclass[11pt]{article}

\usepackage{amsmath} 
\usepackage{amsthm}
\usepackage{amssymb} 
\usepackage{hyperref} 
\usepackage{color} 
\usepackage{enumitem}

\newtheorem{thm}{Theorem}
\newtheorem{inspr}[thm]{}

\newenvironment{definitie}{\begin{itemize}\item[ ]\hspace{-26pt}\bf Definition \rm }{\end{itemize}}
\newenvironment{notatie}{\begin{itemize}\item[ ]\hspace{-26pt}\bf Notation \rm }{\end{itemize}}
\newenvironment{voorbeeld}{\begin{itemize}\item[ ]\hspace{-26pt}\bf Example \rm }{\end{itemize}}
\newenvironment{stelling}{\begin{itemize}\item[ ]\hspace{-26pt}\bf Theorem \rm }{\end{itemize}}
\newenvironment{propositie}{\begin{itemize}\item[ ]\hspace{-26pt}\bf Proposition \rm }{\end{itemize}}
\newenvironment{lemma}{\begin{itemize}\item[ ]\hspace{-26pt}\bf Lemma \rm }{\end{itemize}}
\newenvironment{opmerking}{\begin{itemize}\item[ ]\hspace{-26pt}\bf Remark \rm }{\end{itemize}}
\newenvironment{voorwaarde}{\begin{itemize}\item[ ]\hspace{-26pt}\bf Condition \rm }{\end{itemize}}

\renewcommand{\Bbb}{\mathbb} 

\newcommand{\defin}{\begin{inspr}\begin{definitie}}  %\def already defined
\newcommand{\edefin}{\end{definitie}\end{inspr}}

\newcommand{\notat}{\begin{inspr}\begin{notatie}}  %\not already defined
\newcommand{\enotat}{\end{notatie}\end{inspr}}

\newcommand{\voorb}{\begin{inspr}\begin{voorbeeld}}  %\not already defined
\newcommand{\evoorb}{\end{voorbeeld}\end{inspr}}

\newcommand{\stel}{\begin{inspr}\begin{stelling}}
\newcommand{\estel}{\end{stelling}\end{inspr}}

\newcommand{\prop}{\begin{inspr}\begin{propositie}}
\newcommand{\eprop}{\end{propositie}\end{inspr}}

\newcommand{\lem}{\begin{inspr}\begin{lemma}}
\newcommand{\elem}{\end{lemma}\end{inspr}}

\newcommand{\opm}{\begin{inspr}\begin{opmerking}}
\newcommand{\eopm}{\end{opmerking}\end{inspr}}

\newcommand{\voorw}{\begin{inspr}\begin{voorwaarde}}
\newcommand{\evoorw}{\end{voorwaarde}\end{inspr}}

\newcommand{\bew}{\vspace{-0.3cm}\begin{itemize}\item[ ] \bf Proof\rm: }
\newcommand{\ebew}{\hfill $\qed$ \end{itemize}}

\newcommand{\ssnl}{\vskip 5pt} % Het is noodzakelijk dat er voor de instructie vspace een 
\newcommand{\snl}{\vskip 7pt} % Het is noodzakelijk dat er voor de instructie vspace een lege lijn staat
\newcommand{\nl}{\vskip 12pt} % Kunnen we dit hierin opnemen?

\newcommand{\rood}{\color{red}}
\newcommand{\blauw}{\color{blue}}
\newcommand{\zwart}{\color{black}}

\newcommand{\ot}{\otimes}

\newcommand{\tl}{\triangleleft}
\newcommand{\tr}{\triangleright}

\newcommand{\tussenen}{\qquad\quad\text{and}\qquad\quad}
\newcommand{\tussen}{\qquad\quad\qquad\quad}

\newcommand{\inv}{^{-1}}

\numberwithin{thm}{section}   % Zorgt ervoor dat de nummering bestaat uit ?.?
\numberwithin{equation}{section} % Zorgt ervoor dat de nummering bestaat uit ?.?

\newcommand{\mycomment}[2]{{\blauw {#1}\ssnl}{\zwart{#2}}}
\newcommand{\keepcomment}[1]{}
\newcommand{\oldcomment}[1]{}
\begin{document}

\centerline{\bf \Large Algebraic quantum groups and duality I} 
\vspace{13 pt}
\centerline{\it A.\ Van Daele \rm ($^*$)}
\vspace{20 pt}
{\bf Abstract} 
\nl
Let $(A,\Delta)$ be a finite-dimensional Hopf algebra. The linear dual $B$ of $A$ is again a finite-dimensional Hopf algebra. The \it duality \rm is given by an element $V\in B\ot A$, defined by $\langle V,a\ot b\rangle=\langle a,b\rangle$ where $a\in A$ and $b\in B$. We use $\langle\,\cdot\, , \,\cdot\,\rangle$ for the pairings. In the introduction of this paper, we recall the various  properties of this element $V$ as sitting in the algebra $B\ot A$.
\ssnl
More generally, we can consider \emph{an algebraic quantum group} $(A,\Delta)$. We use the term here for a regular multiplier Hopf algebra with integrals. For $B$ we now take the dual $\widehat A$ of $A$. It is again an algebraic quantum group. In this case, the duality gives rise to an element $V$ in the multiplier algebra $M(B\ot A)$. Still, most of the properties of $V$ in the finite-dimensional case are true in this more general setting.
\ssnl
The \it Heisenberg algebra \rm $C$ is the algebra generated by $A$ and $B$ subject to commutation rules so that $V(a\ot 1)=\Delta(a)V$ holds in the multiplier algebra of $C\ot A$. This algebra acts in a natural way on $A$ and if we consider left multiplication with an element of $A$ in the second factor, we have
\begin{equation*}
V(a\ot a')=\Delta(a)(1\ot a')
\end{equation*}
for all $a,a'\in A$. In other words, $V$ acts as the canonical map $a\ot a'\mapsto \Delta(a)(1\ot a')$  from $A\ot A$ to itself.
\ssnl
The focus in this paper lies further on various aspects of the duality between $A$ and its dual $\widehat A$. Among other things we include a number of formulas relating the  objects associated with an algebraic quantum group and its dual.
\ssnl
This note is meant to give a comprehensive, yet concise  (and sometimes simpler)  account of these known results.
\ssnl
This is part I of a series of three papers on this subject. The case of a multiplier Hopf $^*$-algebra with \it positive integrals \rm is treated in detail in part II \cite{VD-part2} and part III \cite{VD-part3}.
\nl
Date: {\it 26 April 2023} 

\vskip 3cm
\hrule
\vskip 7 pt
\begin{itemize}
\item[($^*$)] Department of Mathematics, KU Leuven, Celestijnenlaan 200B,\newline
B-3001 Heverlee (Belgium). E-mail: \texttt{alfons.vandaele@kuleuven.be}
\end{itemize}

\newpage

\setcounter{section}{-1}  % Dit zorgt ervoor dat we met 0 beginnen voor de inleiding

% Hieronder kunnen de verschillende secties en appendices al of niet opgenomen worden

\section{\hspace{-17pt}. Introduction} \label{sect:introduction}  % \input artikel0.tex%\newpage

First consider two {\it finite-dimensional} Hopf algebras $A$ and $B$ (over the field $\mathbb C$ of complex numbers), together with a non-degenerate pairing $(a,b) \mapsto \langle a,b \rangle$ from the Cartesian product $A\times B$ to $\mathbb C$. We assume that the product in $A$ induces the coproduct on $B$ and vice versa. By this we mean that 
$$\langle\Delta(a),b\ot b'\rangle=\langle a, bb'\rangle
\tussenen
\langle a\ot a',\Delta(b)\rangle= \langle aa',b\rangle$$
for all $a,a'\in A$ and $b,b'\in B$. We use $\Delta$ for the coproduct on $A$ and for the coproduct on $B$. 
Because the pairing is assumed to be non-degenerate, it follows from this that 
$\langle S(a),b\rangle=\langle a,S(b)\rangle$ for all $a,b$ where we denote the antipode with $S$ for $A$ as well as for $B$. Also we have $\langle a,1\rangle=\varepsilon(a)$ and $\langle 1,b\rangle=\varepsilon(b)$ for $a\in A$ and $b\in B$. Again we use $1$ to denote the identity in $A$ and in $B$ and $\varepsilon$ for the counit in the two cases. 
\ssnl
%All of this is well-known and the Hopf algebra $B$
In this case the Hopf algebra $B$ has to be the dual of the Hopf algebra $A$. See e.g.\ Section 7.4 in \cite{R-bk}.
\ssnl
The pairing can be considered as an element $V$ sitting in $B\ot A$ defined by %the formula 
\begin{equation*}
\langle V,a\ot b\rangle=\langle a,b\rangle
\end{equation*}
 for all $a\in A$ and $b\in B$. In the second factor on the left, the original pairing is used. In the first factor, we use the flipped pairing from $B\times A$ to $\mathbb C$.
\ssnl
It is easily shown, using the axioms of a Hopf algebra, that $V$ is invertible in $B\ot A$ and that the inverse is given by $(S\ot \iota)V$ (which is the same as $(\iota\ot S)V$). We use $\iota$ for the identity map, both on $A$ and on $B$. The fact that the product and the coproduct are adjoint to each other is expressed in terms of the element $V$ by the formulas
\begin{equation}(\iota\ot\Delta)V=V_{12}V_{13}
\tussenen
(\Delta\ot\iota)V=V_{13}V_{23}\label{eqn:0.1}
\end{equation}
where we use the standard {\it leg numbering} notation (see further in this introduction).
\ssnl
Next we consider the associated {\it Heisenberg algebra}. It is the algebra generated by $A$ and $B$ subject to the {\it Heisenberg commutation relations}
$$ba=\sum_{(a),(b)} a_{(1)}b_{(2)}\langle a_{(2)},b_{(1)}\rangle$$
for all $a\in A$ and $b\in B$. We use the Sweedler notation for $\Delta(a)$ and for $\Delta(b)$. Again see further in the introduction. We will denote the Heisenberg algebra here by $C$. 
%\ssnl
The linear map from $A\ot B$ to  $C$ given by  $a\ot b\mapsto ab$ is a bijection. 
% and for this reason, we sometimes denote it by $AB$. 
%This justifies the use of $AB$ to denote this algebra. 
\ssnl
In the algebra $C\ot A$ we have the formula
$V(a\ot 1)=\Delta(a)V$
for all $a$. If we combine this formula with the equations in  (\ref{eqn:0.1}), we get the so-called {\it pentagon equation}
\begin{equation}V_{23}V_{12}=V_{12}V_{13}V_{23}.\label{eqn:0.2}
\end{equation}
\nl
The algebra $C$ has a faithful left action on $A$,  given by 
\begin{equation*}
 (ab)\tr x=\sum_{(x)}ax_{(1)}\langle x_{(2)},b\rangle
\end{equation*}
where $a,x\in A$ and $b\in B$. In fact, this can be used to show the injectivity of the map $a\ot b\mapsto ab$ above. Then the element $V$ acts on $A\ot A$ as the canonical map 
\begin{equation*}
T:a\ot a'\mapsto \Delta(a)(1\ot a').
\end{equation*}
With this property, the pentagon equation (\ref{eqn:0.2}) for $V$ follows from the pentagon equation for $T$. The latter combines coassociativity of $\Delta$ with the fact that $\Delta$ is an algebra map.
\nl
All  of these  results mentioned above can be found in the existing literature on Hopf algebras (see further in this introduction). On the other hand, in this paper, we will show that these properties are still valid in a more general situation. The above results can then be obtained as a special case.
\nl
\bf Content of the paper\rm
\nl
The first two sections provide the necessary material needed for the main section, namely  Section \ref{sect:dual}, where we investigate the duality of multiplier Hopf algebras.
\ssnl
In {\it Section} \ref{sect:prel} we recall the objects and relations among these objects for a multiplier Hopf algebra with integrals. In {\it Section} \ref{sect:dual} we consider the dual and we focus on the relations of the objects for the dual and those for the original multiplier Hopf algebra. Some proofs are included. 
\ssnl
We treat the {\it Fourier transform} at the end of this section. It is a distinguished linear map from $A$ to $\widehat A$. 
\ssnl
{\it Section} \ref{sect:dual} contains the main results of this note. We generalize the properties we have above in the finite-dimensional case to the more general case of possibly infinite-dimensional, but  regular multiplier Hopf algebras. In this case, a given multiplier Hopf algebra can be paired with different other multiplier Hopf algebras. 
\ssnl
We will first give notions and results in this more general setting. However, for many of the nicer results, we will assume that the regular multiplier Hopf algebra $A$ has integrals, so that it is an algebraic quantum group, and that $B$ is actually the dual multiplier Hopf algebra $\widehat A$. For the pairing we then take the canonical pairing of $A$ with its dual $\widehat A$. As mentioned already, this generalizes the case of a finite-dimensional Hopf algebra, paired with its dual.
\ssnl
We show that in this situation, the duality $V$ belongs to the multiplier algebra $M(B\ot A)$ and that it satisfies the same formulas as in the much more restrictive finite-dimensional case. The Heisenberg algebra is studied, as well as its action on $A$. Again, this results in the action of $V$ on $A\ot A$ as the canonical map $a\ot a'\mapsto \Delta(a)(1\ot a')$ on $A\ot A$. 
\ssnl
We find how the action of the Heisenberg algebra on $A$ transforms, using the Fourier transform, into an action on $\widehat A$. In particular, we see how the action of the duality, the canonical map, yields another map from $\widehat A\ot \widehat A$ to itself that is of some interest.
% in Section \ref{sect:posint}. 
We also find the formula for the transformation of $V$ under the Fourier transform.
\nl
In the last section, {\it Section} \ref{sect:conclusions},  we draw some conclusions and make further remarks. In particular, we relate this note with the existing literature.
%\mycomment{ Dit mogen we dan niet vergeten. En eventueel de nodige aanpassingen doen. \rood To do!}{} 
 We have to mention that most of the results we formulate and discuss in this paper are known. However, the approach is different in the sense that we start from the pairing and use the Heisenberg algebra as one of the basic tools. This is a more clarifying point of view and therefore in a way also simpler. This note, together with part II \cite{VD-part2} and part III \cite{VD-part3}, is also meant to help readers with a background in algebra to understand the more complicated theory of locally compact quantum groups. 
\nl
\bf Motivation \rm
\nl
While working on a groupoid approach to quantum groups from so-called \emph{matched pairs of closed subgroups} of a locally compact group, we encountered pairings of algebras that, in the most general case, do not fit into the duality framework of algebraic quantum groups. But before studying the necessary generalizations, we found that results where scattered throughout the existing literature. We felt there was a need to give an updated account of these properties before we could proceed and obtain proper and for our purposes, necessary generalizations. We refer here to various papers by M.B.\ Landstad and my self \cite{La-VD3, La-VD4, 
La-VD5} etc. 

\ssnl
This motivated us to write this note. 

\nl
\bf Notations and conventions \rm
\nl
In this note, we only work with algebras over $\Bbb C$. However we believe that many of the results are still true for algebras over other fields.
% This mainly applies to the results in Section 1, as well as for the material of Appendic \ref{sect:appA}, \ref{ sect:appB} and \ref{sect:appC}.
\ssnl
We do not assume that the algebras are unital but we need that the product is \emph{non-degenerate}. In fact, 
the algebras will have local units. Recall that an algebra $A$ is said to have {\it local units} if for any finite set of elements $\{a_1,a_2,\dots,a_n\}$ in $A$, there exists an element $e\in A$ so that $ea_i=a_ie=a_i$ for all $i$. This will imply that the product is non-degenerate. It also follows that the algebras are \emph{idempotent}. This means that any element is a sum of elements $ab$ where both $a$ and $b$ are elements in the algebra. We write this property as $A=A^2$. \keepcomment{Add references!}
\ssnl
If $A$ is a non-degenerate algebra, we use $M(A)$ for the \emph{multiplier algebra} of $A$. When $m$ is in $M(A)$, then by definition we can define $am$ and $mb$ in $A$ for all $a,b\in A$ and we have $(am)b=a(mb)$. The algebra $A$ sits in $M(A)$ as an essential two-sided ideal and $M(A)$ is the largest algebra with identity having this property. If already $A$ has an identity, then $M(A)=A$. If not then $M(A)$ will of course be strictly larger than $A$.
%\begin{commentaar}
%\ssnl
%Misschien moeten we hier een opmerking aan toevoegen over de beperktheid van dit begrip, met een verwijzing naar het idempotent zijn hier?
%\end{commentaar}
\ssnl
A left action $(a,x)\mapsto a\tr x$ of an algebra $A$ on a vector space $X$ is called {\it unital} if any element in $X$ can be written as a sum of elements of the form $a\tr x$ where $a\in A$ and $x\in X$. %We use $\tr$ to denote the action. 
It is called \emph{non-degenerate} if, given $x\in X$, we must have $x=0$ if $a\tr x=0$ for all $a\in A$. If the algebra $A$ has local units, and if the action is unital, for every $x\in X$ there is an element $e\in A$ so that $e\tr x=x$. In particular, in this case the action is automatically non-degenerate. Any unital and non-degenerate action of a non-degenerate algebra extends to a unique unital action of the multiplier algebra (in a canonical way).
\ssnl
If $A$ and $B$ are non-degenerate algebras and if $\alpha:A\to M(B)$ is a homomorphism, it is called {\it non-degenerate} if $\alpha(A)B=B\alpha(A)=B$. In that case, it has a unique extension to a unital homomorphism from $M(A)\to M(B)$. The extension is still denoted by $\alpha$. The same result is true for anti-homomorphisms.
\ssnl
For a non-degenerate algebra $A$, we consider $A\ot A$, the tensor product of $A$ with itself. It is again a non-degenerate algebra and we can consider the multiplier algebra $M(A\ot A)$. The same is true for a multiple tensor product. 
\ssnl
We use $1$ for the identity in any of these  multiplier algebras. On the other hand, we %mostly 
use $\iota$ for the identity map on $A$ (or other spaces).
\ssnl
A linear functional $f$ on $A$ is called {\it faithful} if the bilinear map from $A\times A$ to $\Bbb C$, mapping $(a,b)$ to $f(ab)$, is non-degenerate. So, given $a\in A$, we have that $a=0$ if either $f(ab)=0$ for all $b$ or $f(ba)=0$ for all $b$. 
A faithful  linear functional $f$ is said to have the KMS-property if there is an automorphism $\sigma$ of $A$ such that $f(ab)=f(b\sigma(a))$ for all $a,b$. The automorphism is called the \emph{modular automorphism} or \emph{KMS-automorphism}.
\ssnl
We use $\zeta$ for the flip map on $A\ot A$, as well as for its natural extension to $M(A\ot A)$. 
\ssnl
The {\it leg numbering} notation is used. If e.g.\ $E$ is an element in $M(A\ot A)$, we can consider elements $E_{12}$, $E_{23}$ and $E_{13}$ in $M(A\ot A\ot A)$. For the first two we have $E_{12}=E\ot 1$ and $E_{23}=1\ot E$, whereas for the third one we have $E_{23}=(\iota\ot\zeta)E_{12}$.
\ssnl
For a coproduct $\Delta$ on an algebra $A$, with values in $M(A\ot A)$, 
as we use it in this theory, we assume that $\Delta(a)(1\ot a')$ and $(a\ot 1)\Delta(a')$ are in $A\ot A$ for all $a,a'\in A$. 
\ssnl
We will use the {\it Sweedler notation} for a coproduct on an algebra, also in the case of a multiplier Hopf algebra $(A,\Delta)$ where it is not assumed that $\Delta$  maps $A$ into the tensor product $A\ot A$ but rather in its multiplier algebra $M(A\ot A)$. However, because we assume that $\Delta(a)(1\ot a')$ and $(a\ot 1)\Delta(a')$ are in $A\ot A$ for all $a,a'\in A$, the use of the Sweedler notation is still justified. 
\ssnl 
The use of the Sweedler notation has been introduced in \cite{Dr-VD}. 
We also refer to \cite{VD-tools} and to the more recent paper \cite{VD-sw} (in preparation). 
Crucial for a proper use of the Sweedler notation is to check that the factors are \emph{covered}. 

\ssnl
Remark however that the Sweedler notation is essentially only used as a means to write formulas and equations in a more readable way and not really as a tool to prove results. 

\ssnl
Finally, in order to avoid too many different notations, subscripts, etc., as mentioned before in the finite-dimensional case, we use the same symbol for different objects. 

We hope this will not lead to confusions. In any case, whenever there is some possible doubt, we will be more explicit.
\nl
\bf Basic references \rm
\nl
For the theory of Hopf algebras, we refer to the well-known books by Abe \cite{Ab} and Sweedler \cite{Sw}. See also the more recent work by Radford \cite{R-bk}. The original work on multiplier Hopf algebras is \cite{VD-mha} and for multiplier Hopf algebras with integrals, it is \cite{VD-alg}. The use of the Sweedler notation for multiplier Hopf algebras has been explained in e.g.\ \cite{VD-tools} and more recently in \cite{VD-sw}.
\ssnl
Pairings of multiplier Hopf algebras have been first studied in \cite {Dr-VD}. Actions of multiplier Hopf algebras are studied in \cite{Dr-VD-Z}.
\ssnl
Occasionally, we refer to work on algebraic quantum hypergroups \cite{De-VD2, La-VD4} and papers on weak multiplier Hopf algebras \cite{T-VD-W, VD-W}
\nl
\bf Acknowledgments \rm
\nl
I am very grateful to my coauthor, M.B.\ Landstad and other colleagues and friends, both at the University of Trondheim and the University of Oslo (where part of this note was developed) for the nice and fruitful atmosphere during my regular visits to these 
departments.

 \section{\hspace{-17pt}. Algebraic quantum groups} \label{sect:prel}% \input artikel1.tex

 In this preliminary section, we recall some of the main concepts and its properties encountered in the theory of algebraic quantum groups. We use this term for a multiplier Hopf algebra with integrals as studied first in \cite{VD-mha} and \cite{VD-alg}.
 \ssnl
With any multiplier Hopf algebra $(A,\Delta)$ with integrals are associated several objects. We not only have the coproduct $\Delta$, the counit $\varepsilon$ and the antipode $S$, but also the left and right integrals $\varphi$ and $\psi$ together with  their modular automorphisms $\sigma$ and $\sigma'$. There is the modular element $\delta$ in $M(A)$ relating $\psi$ with $\varphi$ and the scaling constant $\tau$ given by $\varphi\circ S^2=\tau\varphi$.
\ssnl
There are also many relations  between these objects. We collect the most important ones in this  section. We include some proofs of the more special relations.
\ssnl
Furthermore, we have these objects also for the dual $\widehat A$.
Now there are not only the relations of these dual objects among themselves, but moreover also relations between the objects for $A$ and the objects for $\widehat A$. We consider these in Section \ref{sect:dual}.
\ssnl
Throughout we discuss the results and provide references.
\nl
\bf Left and right integrals \rm
\nl
In what follows, we consider a regular multiplier Hopf algebra $(A,\Delta)$ with integrals. We fix a left integral $\varphi$ and we take for the right integral $\psi$ the functional $\varphi\circ S$ where $S$ is the antipode of $(A,\Delta)$. Recall that the integrals are unique, up to a scalar, and that they are faithful linear functionals.
\ssnl
The integrals are invariant in the sense that
\begin{equation}
(\iota\ot\varphi)\Delta(a)=\varphi(a)1
\tussenen
(\psi\ot\iota)\Delta(a)=\psi(a)1 \label{eqn:1.1}
\end{equation}
for all $a$. These equalities hold in $M(A)$. 
\ssnl
The following two formulas are direct consequences of the invariance properties above. We include a proof for completeness.

\prop\label{prop:1.1}
For the left and right integrals, we have
\begin{align}
S((\iota\ot\varphi)(\Delta(a)(1\ot b)))&=(\iota\ot\varphi)((1\ot a)\Delta(b))\label{eqn:1.2}\\
S((\psi\ot\iota)((a\ot 1)\Delta(b)))&=(\psi\ot\iota)(\Delta(a)(b\ot 1))\label{eqn:1.3}
\end{align}
for all $a$ and $b$ in $A$.
\eprop

\bew 
Let $a,b\in A$. Using the basic formulas and the Sweedler notation, we find
\begin{align*}
\sum_{(a)} (S(a_{(1)})\ot 1)\Delta(a_{(2)}b)
&=\sum_{(a)} (S(a_{(1)}) a_{(2)}\ot a_{(3)})\Delta(b)\\
&=\sum_{(a)} (\varepsilon(a_{(1)})1\ot a_{(2)})\Delta(b)\\
&= (1\ot a)\Delta(b).
\end{align*}
We now apply $\iota\ot\varphi$ and use left invariance of $\varphi$ to get
\begin{equation*}
\sum_{(a)} S(a_{(1)})\varphi(a_{(2)}b)
=(\iota\ot\varphi)((1\ot a)\Delta(b)).
\end{equation*}
This proves Equation (\ref{eqn:1.2}). The other equation is proven in a similar way.
\ebew

\ssnl
The formulas (\ref{eqn:1.2}) and (\ref{eqn:1.3}) are important and are used at various places in the paper. 
\ssnl
These formulas are already found in the original paper, see e.g.\ the proof op Proposition 3.11 of \cite{VD-alg}. There the argument is given without using the Sweedler notation but it is essentially the same proof.
\ssnl
As a matter of fact, these formulas are true in more general settings. See e.g.\   
Definition 1.9 and Proposition 2.2 in \cite{De-VD2}, where these formulas are considered for algebraic quantum hypergroups and Proposition 1.5 in \cite{VD-W} where they are encountered in the theory of weak multiplier Hopf algebras. 
\ssnl
The equations in Proposition \ref{prop:1.1} are equivalent with the invariance properties as in Equations (\ref{eqn:1.1}). More precisely, we have the following.

\prop
Assume that $\varphi$ is a linear functional satisfying Equation (\ref{eqn:1.2}) for all $a,b\in A$. Then $\varphi$ is left invariant. Similarly, if $\psi$ is a linear functional satisfying Equation (\ref{eqn:1.3} for all $a,b$, then $\psi$ is right invariant.
\eprop

\bew
If we multiply Equation (\ref{eqn:1.2}) with $S(a')$ in the first factor from the left, we find, for all $a,a',b$ that 
\begin{equation*}
S((\iota\ot\varphi)(\Delta(a)(a'\ot b)))=(\iota\ot\varphi)((S(a')\ot a)\Delta(b)).
\end{equation*}
We know that $A$ has local units (see below for an argument). So we can take for $a$ an element satisfying
\begin{align*}
\Delta(a)(a'\ot b))&=a'\ot b\\
(S(a')\ot a)\Delta(b)&=(S(a')\ot 1)\Delta(b).
\end{align*}
We use that $A\ot A\subseteq \Delta(A)(1\ot A)$ and that $(A\ot 1)\Delta(A)\subseteq A\ot A$.
Then we get $S(a')\varphi(b)=S(a')((\iota\ot\varphi)\Delta(b))$ and left invariance of $\varphi$ follows.
\ssnl
In a similar way, we can prove the result for right invariant functionals.
\ebew

In a recent paper on finite quantum hypergroups (\cite{La-VD4} there is a discussion about the relation of the formulas (\ref{eqn:1.2}) and (\ref{eqn:1.3}) with the invariance properties stated in (\ref{eqn:1.1}). On the other hand, these arguments will only work when $A$ has an identity. 

\keepcomment{See remark in the tex file for future adaptations.}
\ssnl
The proofs are of course easier when $A$ has an identity. Then we simply can take for $a$ the identity in Equation (\ref{eqn:1.3}) while in Equation (\ref{eqn:1.5}) we take for $b$ the identity. In general, we need local units. They exist for any multiplier Hopf algebras as we see below.

\prop\label{prop:1.3} 
Let $(A,\Delta)$ be any multiplier Hopf algebra. Then $A$ has local units.
\eprop

\bew
i) Let $a\in A$. We claim that $a\in Aa$. To prove this, assume that $\omega$ is a linear functional on $A$ such that $\omega(ba)=0$ for all $b\in A$. Then
\begin{equation*}
(\iota \ot\omega)((c\ot d)\Delta(b)(1\ot a))=0
\end{equation*}
for all $b,c,d$ because $(c\ot d)\Delta(b)\in A\ot A$. Because already $(1\ot d)\Delta(b)(1\ot a)\in A\ot A$ and because the product in $A$ is non-degenerate, we can cancel $c$ and obtain that also
\begin{equation*}
(\iota \ot\omega)((1\ot d)\Delta(b)(1\ot a))=0
\end{equation*}
for all $b,d$. Write $\Delta(b)(1\ot a)=\sum_ip_i\ot q_i$, replace $d$ by $S(p_i)$ and take sums to get
\begin{equation*}
\sum_i\omega (S(p_i)q_i)=0.
\end{equation*}
Because $\sum_i S(p_i)q_i=\varepsilon(b)a$, with $\varepsilon(b)=1$ we find $\omega(a)=0$. This shows that $a\in Aa$ for all $a$.
\ssnl
ii) In a similar way we get $a\in aA$ for all $a$. 
\ssnl
iii) We see that for each $a$ there exist elements $e,f$ in $A$ satisfying $a=ea=af$. It can be shown that then local units exist (see e.g.\ \cite{Ve}). 
\ebew

Here are some comments on the history of this result. \keepcomment{We have this result also in the paper with Landstad on polynomial functions. If this one is put on the arxiv first, we can remove it there and refer to here.}

\opm
The first result of this type is found in \cite{Dr-VD-Z}. In Proposition 2.2 of that paper, it is shown that any \emph{regular} multiplier Hopf algebra has left and right local units. In Proposition 2.6 of the same paper, it is proven that any regular multiplier Hopf algebra with integrals has two-sided local units. When that paper was published (1999), it was still open if these properties would still hold for any multiplier Hopf algebra.
\ssnl
In \cite{VD-Z},  it is shown in Proposition 1.2 of that paper that in any multiplier Hopf algebra $(A,\Delta)$, given an element $a\in A$, there are elements $e,f$ in $A$ satisfying $ea=af=a$. The proof is as in the first part of the proof above of Proposition \ref{prop:1.3}. This paper was written later than \cite{Dr-VD-Z} but appeared also in 1999. At that time, the result of \cite{Ve} was not yet available and so it was still open whether or not any multiplier Hopf algebra had local units.
\ssnl
Later, the existence of local units has also been obtained for algebraic quantum hypergroups and weak multiplier Hopf algebras. In the first case, the property is found in Proposition 1.6 of \cite{De-VD} while for the second case, it is proven in Proposition 2.14 of \cite{VD-W}.
\eopm

The Equations\ref{eqn:1.2}  can be used to prove the following property of left integrals.

\prop\label{prop:1.5a}
A left integral is faithful. It is also unique up to a scalar.
\eprop

\bew
We sketch the proof. 
\ssnl
i) For the first statement, we assume that $a\in A$ and that $\varphi(ab)=0$ for all $b$. Then using Equation \ref{eqn:1.2} we get that $\sum_{(a)}a_{(1)}\varphi(a_{(2)}b)=0$ for all $b$. We now apply the $\Delta$, the inverse of the antipode and replace $b$ by $S\inv(a_{(2})b$ to get 
\begin{equation*}
\sum_{(a)}a_{(1)}\varphi(a_{(3)}S\inv(a_{(2)})b)=a\varphi(b)=0.
\end{equation*}
This implies that $a=0$. In a similar way we get $a=0$ if $\varphi(ba)=0$ for all $b$. 
\ssnl
ii) To prove uniqueness, let $\varphi'$ be any other left invariant functional. Apply $\varphi'$ to Equation \ref{eqn:1.2} and use that $\varphi'\circ S$ is right invariant. Then we get
\begin{equation*}
\varphi'(S(a))\varphi(b)=\varphi(ac)
\end{equation*}
where $c=(\varphi'\ot\iota)\Delta(b)$. From the faithfulness of $\varphi$ we can conclude that there is a multiplier $\delta$ of $A$ so that $(\varphi'\ot\iota)\Delta(b)=\varphi(b)\delta$. Then we get $\varphi'(S(a))=\varphi(a\delta).$ From this uniques follows.
\ebew

We see that along the way, we also have shown the existence of the modular element $\delta$.

%{\blauw 
% We may have  refer to one of the papers with Magnus. Also refer to the paper on wmha3.
% If we finish that paper}
\nl
\bf The scaling constant and the modular element \rm
\nl
An immediate consequence of the uniqueness of the integrals is the following.
\prop\label{prop:1.4}
There is a non-zero complex number $\tau$ satisfying $\varphi\circ S^2=\tau\varphi$. We also have $\psi\circ S^2=\tau\psi$.
\eprop

The number $\tau$ is called the {\it scaling constant}. There are examples where $\tau\neq 1$.  See Proposition 5.11 in \cite{VD-alg}. 
% {\blauw Discuss further and give references.} 
%see the item `Notes and references' at the end of this appendix. {\blauw To do!}
\nl
Next we have the relation of the left integral with the right integral.
\prop\label{prop:1.5}
There is a unique invertible element $\delta$ in the multiplier algebra $M(A)$ of $A$ satisfying and characterized by 
\begin{equation}
\varphi(S(a))=\varphi(a\delta) \tussenen \varphi(S\inv(a))=\varphi(\delta a) \label {eqn:1.4 }
\end{equation}
for all $a\in A$. It satisfies
\begin{equation}
(\varphi\ot\iota)\Delta(a)=\varphi(a)\delta
\tussenen
(\iota\ot\psi)\Delta(a)=\psi(a)\delta^{-1}\label{eqn:1.5}
\end{equation} 
for all $a$. Moreover $\Delta(\delta)=\delta\ot\delta$, $\varepsilon(\delta)=1$ and $S(\delta)=\delta\inv$.
%i.e.\ $\delta$ is group-like. 
\eprop

The existence of $\delta$ satisfying $\varphi(S(a))=\varphi(a\delta)$ has been obtained in the proof of Proposition \ref{prop:1.5a}.
\ssnl
We can add some more properties. If we let $\psi=\varphi\circ S$ we have, for all $a$,  
\begin{equation}
\psi(S(a))=\varphi(S^2(a))=\tau\varphi(a)=\tau\psi(a\delta\inv)\label{eqn:1.6}
\end{equation}.

We see that $\delta$ is a group-like element in the multiplier algebra $M(A)$. In the algebraic literature on the subject, when the algebra is finite-dimensional, this element is often called the \emph{distinguished group-like element}, see e.g.\ \cite{Ra}. \oldcomment{Referentie verifieren}{}
\ssnl
The two equalities in Equation \ref{eqn:1.4 } are equivalent with each other. Indeed, let $b\in A$ and  $a=S\inv(b)\delta\inv$ then $S(a)=\delta b$ and from $\varphi(S(a))=\varphi(a\delta)$ we get $\varphi(\delta b)=\varphi(S\inv(b))$. Similarly in the other direction.
\ssnl
%\blauw Discuss this also in `Notes and references'.
%\ssnl}

For the last equation, we use the unique extension of the coproduct $\Delta$ to the multiplier algebra. The formulas (\ref{eqn:1.5}) hold in the multiplier algebra $M(A)$.  This invertible multiplier $\delta$ is called the {\it modular element}. From the fact that $\delta$ is group-like, we get $\varepsilon(\delta)=1$ and $S(\delta)=\delta^{-1}$. Again we use the extensions of the counit and the antipode to the multiplier algebra. In particular $S^2(\delta)=\delta$.
\nl
\bf The modular automorphisms \rm
\nl
Finally, we have the modular properties of the integrals.
\prop\label{prop:1.6}
There exist unique automorphisms $\sigma$ and $\sigma'$ of $A$ satisfying and characterized by
\begin{equation*}
\varphi(ab)=\varphi(b\sigma(a))
\tussenen
\psi(ab)=\psi(b\sigma'(a))
\end{equation*}
for all $a,b\in A$. We have $\varphi(\sigma(a))=\varphi(a)$ and $\psi(\sigma'(a))=\psi(a)$ for all $a$.
\eprop

These automorphisms are called the \emph{modular automorphisms}. 
\ssnl
The existence of $\sigma$ and $\sigma'$ as linear maps is not obvious. The result follows from a systematic use of the basic formulas (\ref{eqn:1.2}) and (\ref{eqn:1.3} in Proposition \ref{prop:1.1}.  For such a treatment, we refer to \cite{VD-W} where this result is proven for weak multiplier Hopf algebras with integrals. See Proposition 1.7 in \cite{VD-W}.
\ssnl
It follows easily from the faithfulness that these maps $\sigma$ and $\sigma'$ are homomorphisms. In fact, they are automorphisms.
Injectivity also follows from the faithfulness of the integrals while surjectivity is again non-trivial.  To prove the invariance of the integrals under the modular automorphisms, observe e.g.\  that
\begin{equation*}
\varphi(ab)=\varphi(b\sigma(a))=\varphi(\sigma(a)\sigma(b))=\varphi(\sigma(ab))
\end{equation*}
for all $a,b$. As $A$ is idempotent, we get $\varphi(\sigma(a))=\varphi(a)$ for all $a\in A$. A similar argument is used for $\psi$.

\opm
A finite-dimensional unital algebra with a faithful linear functional is called a Frobenius algebra. In that case, the inverse of the modular automorphism is called the Nakayama automorphism (see \cite{Nak}). The terminology we use comes from the theory of locallly compact groups and operator algebras (see e.g.\ \cite{Pe}). Remark that if the algebra is finite-dimensional and has a faithful functional, then the modular automorphism always exists. This is not true for infinite-dimensional algebras. % Later referentie naar werk met Joost opnemen.
\eopm

We now formulate a number of relations involving these modular automorphisms.
\ssnl
First we have the following relations between the two modular automorphisms, the antipode and the modular element.

\prop\label{prop:1.8}
We have $\sigma(S(a))=S({\sigma'}^{-1}(a))$ and $\sigma(S^{-1}(a))=S^{-1}({\sigma'}^{-1}(a))$ for all $a$. Also $\sigma'(a)=\delta\sigma(a)\delta^{-1}$ for all $a$. Finally we find 
\begin{equation*}
\sigma(\delta)=\tau^{-1}\delta
\tussenen
\sigma'(\delta)=\tau^{-1}\delta.
\end{equation*}
\eprop

The first two equations are immediate consequences of the formulas $\psi(a)=\varphi(S(a))$ and $\varphi(a)=\tau\psi(S(a))$ for all $a$. On the other hand, the equality $\psi(a)=\varphi(a\delta)$ will give the last formula.
\ssnl
We use these equalities to prove the following.

\prop\label{prop:1.9}
The automorphisms $S^2$, $\sigma$ and $\sigma'$ mutually commute with each other.
\eprop

\bew
For any $a$ we have
\begin{equation*}
\sigma(S^2(a))=S({\sigma'}^{-1}(S(a)))=\sigma(S^2(a)).
\end{equation*}
For the first equality we have used the first equality in the previous proposition while for the second one, we used the second equality in that proposition. So we obtain that $\sigma$ and $S^2$ commute. 
\ssnl
Using that $S^2(\delta)=\delta$ and the relation $\sigma'(a)=\delta\sigma(a)\delta^{-1}$ we get that also $S^2$ and $\sigma'$ will commute. Finally, using that $\sigma(\delta)=\tau^{-1}\delta$, again from this equation, we will find that $\sigma$ and $\sigma'$ commute. 
\ebew

In the following proposition, we look for the behavior of the counit under the modular automorphisms.

\prop\label{prop:1.10} 
For all $a$ in $A$ we have
\begin{equation}
\varepsilon(\sigma(a))=\varepsilon(\sigma'(a))
\tussenen
\varepsilon(\sigma^{-1}(a))=\varepsilon({\sigma'}^{-1}(a)) \label{eqn:1.7}.
\end{equation}
\eprop

\bew
That $\varepsilon\circ\sigma=\varepsilon\circ\sigma'$ follows from $\sigma'(a)=\delta\sigma(a)\delta^{-1}$ by applying $\varepsilon$ and using that $\varepsilon(\delta)=1$. In a similar way we get $\varepsilon(\sigma^{-1}(a))=\varepsilon({\sigma'}^{-1}(a))$.
\ebew

Finally we have the relations of the modular automorphisms with the coproduct.

\prop\label{prop:1.11} 
We have 
\begin{align}
\Delta(\sigma(a))&=(S^2\ot\sigma)\Delta(a) \label{eqn:1.8}\\
\Delta(\sigma'(a))&=(\sigma'\ot S^{-2})\Delta(a)\label{eqn:1.9}\\
\Delta(S^2(a))&=(\sigma\ot{\sigma'}^{-1})\Delta(a)\label{eqn:1.10}
\end{align}
for all $a$.
\eprop

\bew
i) The formula in Equation (\ref{eqn:1.8}) is obtained by applying Equation (\ref{eqn:1.2}) twice while Equation (\ref{eqn:1.9}) follows in a similar way from Equation (\ref{eqn:1.3}. 
\ssnl
ii) If we apply $\varepsilon$ on the second leg of Equation (\ref{eqn:1.8}) and on the first leg of Equation (\ref{eqn:1.9}) we find
\begin{equation}
(\iota\ot\varepsilon\circ\sigma)\Delta(a)=S^{-2}\sigma(a)
\tussenen
(\varepsilon\circ\sigma'\ot\iota)\Delta(a)=S^2\sigma'(a)\label{eqn:1.11}
\end{equation}
for all $a$. From Equation (\ref{eqn:1.7})  we find
\begin{equation*}
(\iota\ot\varepsilon\circ\sigma\ot\iota)\Delta^{(2)}(a)=(\iota\ot\varepsilon\circ\sigma'\ot\iota)\Delta^{(2)}(a)
\end{equation*}
and if we combine this with the equations in (\ref{eqn:1.11}) and use coassociativity of $\Delta$ we find
\begin{equation*}
(S^{-2}\sigma\ot\iota)\Delta(a)=(\iota\ot S^2\sigma')\Delta(a)
\end{equation*}
and so
\begin{equation*}
(\sigma\ot{\sigma'}^{-1})\Delta(a)=(S^2\ot S^2)\Delta(a)=\Delta(S^2(a))
\end{equation*}
for all $a$. This proves the formula in (\ref{eqn:1.10}).
\ebew

All these formulas, except Equation (\ref{eqn:1.10}), are already present in the original paper on multiplier Hopf algebras with integrals \cite{VD-alg}. The formula (\ref{eqn:1.10}) was first obtained in Lemma 3.10 of \cite{Ku-VD} in the case of a multiplier Hopf $^*$-algebra with positive integrals. However later  a more direct proof has been given, valid in the more general case (see e.g.\ Proposition A.6 in the appendix of \cite{La-VD2} and Proposition 2.7 in \cite{De-VD2}). 
\nl

\oldcomment{
The following still has to be completely checked:
\begin{itemize}\setlength\itemsep{0pt}
\item References for these statements - To be completed here and there
\item Existence of local units, $A^2=A$ - More or less ok
\item About the terminology - More or less ok
\item modular element is distinguished group like element - More or less ok
\item modular automorphisms are Nakayama automorphisms - More or less ok
\item existence of the modular automorphisms - More or less ok
\end{itemize}
}{}
%$\tau$ can be non-trivial

 \section{\hspace{-17pt}. The dual of an algebraic quantum group} \label{sect:dual}% \input artikel2.tex

As in the previous section, we  consider a multiplier Hopf algebra $(A,\Delta)$ with a left integral $\varphi$. For the right integral $\psi$ we take $\varphi\circ S$. The dual $\widehat A$ is defined as a subspace of the dual space of $A$. They are elements of the form $\varphi(\,\cdot\,a)$ where $a\in A$. In what follows, we use $B$ for $\widehat A$ and consider the obvious pairing of $A$ with $B$ given by evaluation of the elements of $B$ in $A$.
\ssnl
 The coproduct on $A$ gives a product on $B$ while the product on $A$ induces a coproduct on $B$. This coproduct is again denoted by $\Delta$. As a matter of fact, we have
 \begin{equation*}
\langle a,bb'\rangle=\langle \Delta(a),b\ot b'\rangle\tussen\langle a\ot a',\Delta(b)\rangle=\langle aa',b\rangle
\end{equation*}
for $a,a'\in A$ and $b,b'\in B$. The formulas make sense because it turns out to be possible to extend the pairing of $A\ot A$ with $B\ot B$ to $M(A\ot A)\times (B\ot B)$ and to $(A\ot A)\times M(B\ot B)$. 
\ssnl
All these results require some work. We will consider these extensions and these formulas in the more general setting of a pairing between regular multiplier Hopf algebras in the first item of the next section. It is more natural to do this for such a more general pairing.
 \nl
 The pair $(B,\Delta)$ is also an algebraic quantum group. We use  again $S$ for the antipode and $\varepsilon$ for the counit on the dual. For the other objects of the dual, we use the same symbols, covered with a hat. For the construction and the properties of the dual, we refer to \cite{VD-alg}.
 \ssnl
 
 \nl
 \bf The dual objects \rm
 \nl
 First we have the integrals on the dual.
 
 \prop\label{prop:2.1}
Define $\widehat\psi$ on $B$ by $\widehat \psi(b)=\varepsilon(a)$ if $b=\varphi(\,\cdot\,a)$. This defines a right integral on $B$. On the other hand, if we define $\widehat\varphi$ on $B$ by $\widehat\varphi(b)=\varepsilon(a)$ if $b=\psi(S(\,\cdot\,)a)$, then $\widehat\varphi$ is a left integral on $B$.
\eprop

\bew
i) Take $a\in A$ and let $b=\varphi(\,\cdot\,a)$. For $a',a''$ in $A$ we have
\begin{equation*}
\langle a'\ot a'',\Delta(b) \rangle=\langle a'a'',b\rangle=\varphi(a'a''a).
\end{equation*}
It follows that 
\begin{equation*}
(\iota\ot\langle a'',\,\cdot\, \rangle)\Delta(b)=\varphi(\,\cdot\,a''a).
\end{equation*}
Then, from the definition of $\widehat\psi$,  
\begin{align*}
\widehat\psi((\iota\ot\langle a'',\,\cdot\, \rangle)\Delta(b))
&=\varepsilon(a''a)=\varepsilon(a'')\varepsilon(a)\\
&=\varepsilon(a'')\widehat\psi(b)=\widehat\psi(b)\langle a'', 1\rangle.
\end{align*}
This proves that $(\widehat\psi\ot\iota)\Delta(b)=\widehat\psi(b) 1$.
\vskip 3pt
ii) Again take $a\in A$ and let $b=\psi(S(\,\cdot\,)a)$. For $a',a''$ in $A$ we now find
\begin{equation*}
\langle a'\ot a'',\Delta(b) \rangle=\langle a'a'',b\rangle=\psi(S(a'a'')a)=\psi(S(a'')S(a')a).
\end{equation*}
Then 
\begin{equation*}
(\langle a',\,\cdot\,\rangle\ot\iota)\Delta(b)=\psi(S(\,\cdot\,)S(a')a).
\end{equation*}
It follows that
\begin{align*}
\widehat\varphi((\langle a',\,\cdot\,\rangle\ot\iota)\Delta(b))
&=\varepsilon(S(a')a)=\varepsilon(S(a'))\varepsilon(a)\\
&=\varepsilon(a')\varepsilon(a)=\widehat\varphi(b)\langle a',1 \rangle.
\end{align*}
This proves that $(\iota\ot\widehat\varphi)\Delta(b)=\widehat\varphi(b)1$.
\ebew

This result is found in Proposition 4.8 of \cite{VD-alg}. With these choices we get the following relation.

\prop\label{prop:2.2} 
Let $\varphi$ be left integral on $A$. Define $\psi$ on $A$ by $\psi(a)=\varphi(S(a))=\varphi(a\delta)$. Associate $\widehat\psi$ and $\widehat\varphi$ as in the previous proposition. Then $\widehat\psi(b)=\widehat\varphi(S(b))$ for all $b\in B$.
\eprop

\bew
Let $b=\varphi(\,\cdot\,c)$ for $c\in A$. Then, as $\psi=\varphi(\,\cdot\,\delta)$, we get 
\begin{equation*}
\langle a,S(b)\rangle=\langle S(a),b\rangle=\varphi(S(a)c)=\psi(S(a)c\delta\inv).
\end{equation*}
It follows that $S(b)=\psi(S(\,\cdot\,)c\delta\inv)$ so that 
\begin{equation*}
\widehat\varphi(S(b))=\varepsilon(c\delta\inv)=\varepsilon(c)=\widehat\psi(b).
\end{equation*}
We have used that $\varepsilon(\delta)=1$, see Proposition \ref{prop:1.5}.
\ebew

We could have used this result to prove the second part of Proposition \ref{prop:2.1}.
\ssnl
In what follows we keep the assumption $\psi=\varphi\circ S$ on $A$ and we take for  $\widehat\varphi$ and $\widehat\psi$  the associated left and right integrals on $B$ as in Proposition \ref{prop:2.1}. It is important to notice that also for the dual, we have the relation  $\widehat\psi=\widehat\varphi\circ S$. These choices are therefore compatible with each other.
 
\nl
\bf Relations among objects of $(A,\Delta)$ and its dual\rm
\nl
In the first place, we have the relations of the dual objects associated with the dual  among themselves as we have them for the objects of the original pair $(A,\Delta)$.  But now, there are also a variety of relations of the objects of $(A,\Delta)$ with the dual objects. These results are not contained in the original paper \cite{VD-alg}, but they are found e.g., in fact in greater generality, in Section 4 of \cite{De-VD2}.
\ssnl
For the following set of formulas, we need the extension of the pairing of $A$ with $B$ to $A\times M(B)$. As mentioned already, we will treat these extended pairing in the next section. They satisfy
\begin{equation*}
\langle a,bm\rangle=\langle a\tl b,m \rangle
\tussenen
\langle a,mb\rangle=\langle b\tr a,m \rangle
\end{equation*}
where 
\begin{equation*}
a\tl b=(\langle \,\cdot\, , b\rangle\ot\iota) \Delta(a)
\tussenen
b\tr a=(\iota\ot \langle \,\cdot\, ,b\rangle) \Delta(a)
\end{equation*}
for $a\in A$, $b\in B$ and $m\in M(B)$. In the case we treat here, so when $B$ is equal to $\widehat A$, these expressions are well-defined because $b$ has the form $\varphi(\,\cdot\, c)$ for some $c\in A$.
\oldcomment{\rood Include a reference. This is treated in the next section in greater generality. See Remark \ref{opm:3.3}. We have to include this earlier as we are using it before. But it is ok to treat it in the next section. This is more general.}{}

\prop\label{prop:2.3}
For all $a\in A$ and $b\in B$ we have
\begin{align}
\langle a, \widehat\delta\rangle &=\varepsilon(\sigma^{-1}(a)) =\varepsilon({\sigma'}^{-1}(a))\label{eqn:2.1} \\
\langle a, \widehat\delta^{-1} \rangle &=\varepsilon(\sigma(a))=\varepsilon(\sigma'(a)).\label{eqn:2.2}
\end{align}
In these two equations we use the extension of the pairing to $A\times M(B)$.
\eprop

\bew
i) Let $a,a'\in A$ and $b,b'\in B$ and assume that $b=\varphi(\,\cdot\,c)$ for $c\in A$. Because the coproducts are adjoints of the products we get 
$$\langle a\ot a',(b'\ot 1)\Delta(b) \rangle=\langle \Delta(a)(1\ot a'),b'\ot b\rangle.$$
 Then, using the Sweedler notation,  we have
\begin{align*}
\langle a\ot a',(b'\ot 1)\Delta(b) \rangle
% &=\langle \Delta(a)(1\ot a'),b'\ot b\rangle\\% \label{eqn:2.1}\\
&=\sum_{(a)} \langle a_{(1)},b'\rangle\varphi(a_{(2)}a'c)\\%\label{eqn:2.4}\\
&=\sum_{(a)} \langle a_{(1)},b'\rangle\varphi( a'c\sigma(a_{(2)})).%\label{eqn:2.2}
\end{align*}
It follows that 
\begin{equation*}
\langle a\ot \,\cdot\, ,(b'\ot 1)\Delta(b) \rangle=\sum_{(a)} \langle a_{(1)},b'\rangle\varphi( \,\cdot\,c\sigma(a_{(2)})).
\end{equation*}
Now apply $\widehat\psi$. For the left hand side, using the second of the Equations (\ref{eqn:1.5}) of Proposition \ref{prop:1.5}, we find
$\widehat\psi(b)\langle a,b'\widehat \delta\inv\rangle$ while for the right hand side we get
\begin{align*}
\sum_{(a)} \langle a_{(1)},b'\rangle\varepsilon(c\sigma(a_{(2)}))
&=\sum_{(a)} \langle a_{(1)},b'\rangle\varepsilon(c)\varepsilon(\sigma(a_{(2)}))\\
&=\sum_{(a)} \langle a_{(1)},b'\rangle\widehat\psi(b)\varepsilon(\sigma(a_{(2)})).
\end{align*}
We find that 
\begin{equation*}
\langle a,b'\widehat\delta\inv\rangle=\sum_{(a)} \langle a_{(1)},b'\rangle\varepsilon(\sigma(a_{(2)})).
\end{equation*}
From the definition of the extension of the pairing from $A\times B$ to $A\times M(B)$ we see that
$\langle a,\widehat\delta\inv\rangle=\varepsilon(\sigma(a))$. This proves the first equality of (\ref{eqn:2.1}). For the second one, we then use the result of Proposition \ref{prop:1.10}.
 \ssnl
ii) To prove (\ref{eqn:2.2}) we can use that $S(\delta)=\delta\inv$. Then we get
 \begin{equation*}
\langle a,\widehat\delta\rangle=\langle S\inv(a),\widehat\delta\inv\rangle
=\varepsilon(\sigma(S\inv(a)))=\varepsilon(S\inv({\sigma'}\inv(a)))=\varepsilon({\sigma'}\inv(a)).
\end{equation*}
Here we have used Proposition \ref{prop:1.8}.
\ebew

To prove the next result, we need the following formula. It is like Plancherel's formula. We will explain this  in the item on the Fourier transform, later in this section.  

\prop\label{prop:2.4}
If $b=\varphi(\,\cdot\,c)$ and $b'=\varphi(\,\cdot\,c')$ then $\widehat\psi(bb')=\varphi(S\inv (c')c)$.
\eprop
\bew
We have, using Equation (\ref{eqn:1.2}) of Proposition \ref{prop:1.1},
\begin{align*}
\langle a,bb'\rangle
&=\sum_{(a)}\langle a_{(1)},b\rangle\varphi(a_{(2)}c') \\
&=\sum_{(c')}\langle S\inv(c'_{(1)}),b\rangle\varphi(ac'_{(2)})
\end{align*}
and so
\begin{equation*}
\widehat\psi(bb')=\sum_{(c')}\langle S\inv(c'_{(1)}),b\rangle\varepsilon(c'_{(2)})=\langle S\inv(c'),b\rangle=\varphi(S\inv(c')c).
\end{equation*}
\ebew

From this result, combined with the formulas in Proposition \ref{prop:2.3}, we can obtain the following.

\prop
For all $a\in A$ and $b\in B$ we have
\begin{align}
\langle a, \widehat\sigma (b)\rangle &= \langle S^2(a)\delta^{-1}, b\rangle \label{eqn:2.3}\\
\langle a, \widehat\sigma'(b)\rangle &= \langle \delta^{-1}S^{-2}(a), b\rangle\label{eqn:2.4}
\end{align}
\eprop
\bew
i) Let $b=\varphi(\,\cdot\,c)$, $b'=\varphi(\,\cdot\,c')$ and $\widehat\sigma'(b)=\varphi(\,\cdot\,c'')$. From $\widehat\psi(bb')=\widehat\psi(b'\widehat\sigma'(b))$, using Proposition \ref{prop:2.4}, we get

\begin{align*}
\varphi(S\inv(c')c)
&=\varphi(S\inv(c'')c')\\
&=\varphi(S\inv(S(c')c''))\\
&=\varphi(\delta S(c')c'').
\end{align*}
It follows that 
%\begin{equation*}
$\langle S\inv(c'),b\rangle=\langle\delta S(c'),\widehat\sigma'(b)\rangle$.
%\end{equation*}
With $a=\delta S(c')$ we find 
$$S\inv(c')=S\inv(S\inv(\delta\inv a))=\delta\inv S^{-2}(a)$$
 and 
%\begin{equation*}
$\langle a,\widehat\sigma'(b)\rangle= \langle\delta\inv S^{-2}(a)\rangle$.
%\end{equation*}
This proves Equation (\ref{eqn:2.4}).
\ssnl
ii) To prove Equation (\ref{eqn:2.3}), we can use a similar argument. However we will use a different method to illustrate some other relations.
\ssnl
First we claim that $\langle a,{\widehat{\sigma'}}\inv(b)\rangle=\langle \delta S^2(a), b\rangle$. This simply follows from 
\begin{equation*}
\langle a, b \rangle=\langle a, \widehat{\sigma'}({\widehat{\sigma'}}\inv(b))\rangle=\langle \delta\inv S^{-2}(a), {\widehat{\sigma'}}\inv(b)\rangle
\end{equation*} by replacing $\delta\inv S^{-2}(a)$ by $c$ so that $a=\delta S^2(c)$ and $\langle c,{\widehat{\sigma'}}\inv(b)\rangle=\langle \delta S^2(c), b\rangle$. %This proves the claim.
\ssnl
Next we use that $S\inv(\widehat\sigma(b))={\widehat{\sigma'}}\inv(S\inv(b))$. Then we get
\begin{align*}
\langle a,\widehat\sigma(b)\rangle
&=\langle S(a),S\inv(\widehat\sigma(b))\rangle\\
&=\langle S(a),{\widehat{\sigma'}}\inv(S\inv(b)) \rangle \\
&=\langle \delta S^3(a),S\inv(b)\rangle \\
&=\langle  S^2(a)\delta\inv,b\rangle.
\end{align*}
This proves the Equation (\ref{eqn:2.3}).
\ebew

Another interesting relation is Radford's formula. 
If $\delta$ and $\widehat \delta$ denote the modular elements in $M(A)$ and $M(\widehat A)$ respectively, then
$$S^4(a)=\delta^{-1}(\widehat \delta\tr a\tl \widehat\delta^{-1})\delta$$
for all $a\in A$. Here $\tr$ stands for the left action of $B$ on $A$ induced from the pairing while $\tl$ is the right action. Recall that these actions commute. Both are extended to the multiplier algebra.
\ssnl
Some of these formulas can be found already in the original papers on the subject (see e.g.\ \cite{VD-mha} and \cite{VD-alg}. In \cite{De-VD-W} more formulas are found an are treated in a systematic way. 
\nl
\bf The Fourier transform \rm
\nl
We now formulate  two results about the Fourier transform. Remember, given the left integral $\varphi$ on $A$,  the choices made for the right integral $\psi$ on $A$ and the integrals $\widehat\varphi$ and $\widehat\psi$ on $B$, see Proposition \ref{prop:2.2} and the remark following it. %\mycomment{Choice of the integrals}

\prop\label{prop:2.6}
The inverse of the map $a\mapsto \varphi (\,\cdot\,a)$ from $A$ to $B$ is the map $b\mapsto \widehat\psi( S(\,\cdot\,)b)$ from $B$ to $A$.
\eprop

\bew
Take $a'\in A$ and $b'\in B$. Then
\begin{align*}
\langle a',S(b')b\rangle
&=\sum_{(a')} \langle a'_{(1)}, S(b')\rangle\langle a'_{(2)},b\rangle \\
&=\sum_{(a')} \langle S(a'_{(1)}),b'\rangle\varphi(a'_{(2)}a) \\
&=\sum_{(a)} \langle a_{(1)},b'\rangle\varphi(a'a_{(2)}). 
\end{align*}
For the last equality, we have used Equation (\ref{eqn:1.2})
\begin{equation*}
S((\iota\ot\varphi)(\Delta(a')(1\ot a)))=(\iota\ot\varphi)((1\ot a')\Delta(a)).
\end{equation*}
We see that $S(b')b=\sum_{(a)} \langle a_{(1)},b'\rangle\varphi(\,\cdot\, a_{(2)})$. By the definition of $\widehat \psi$ we find
\begin{equation*}
\widehat \psi(S(b')b)=\sum_{(a)} \langle a_{(1)},b'\rangle \varepsilon(a_{(2)})=\langle a,b'\rangle.
\end{equation*}
It follows that $\widehat \psi(S(\,\cdot\,)b)=a$ if $b=\varphi(\,\cdot\, a)$ for $a\in A$. This proves the result.
\ebew

This result is found in \cite{VD-alg} where it is used to prove that the dual of $\widehat A$ is canonically isomorphic with $A$, see Theorem 4.12 in \cite{VD-alg}.
\ssnl
This is the first option. Here is the other one.

\prop\label{prop:2.7}
The inverse of the map $a\mapsto \psi(S(\,\cdot\,)a)$ from $A$ to $B$ is the map \newline $b\mapsto \widehat\varphi(\,\cdot\,b)$.
\eprop

\bew
Let $a\in A$ and $b=\psi(S(\,\cdot\,)a)$. For all $d\in B$ and  $x\in A$ we have
\begin{align*}
\langle x,db\rangle
&= \sum_{(x)} \langle x_{(1)},d\rangle \langle x_{(2)},b\rangle \\
&= \sum_{(x)} \langle x_{(1)},d\rangle \psi(S(x_{(2)})a) \\
&= \sum_{(S(x))} \langle S(x)_{(2)},S\inv(d)\rangle \psi((S(x))_{(1)}a) \\
&= \sum_{(a)} \langle S(a_{(2)}),S\inv(d)\rangle \psi(S(x)a_{(1)}).
\end{align*}
We have used Equation (\ref{eqn:1.3}) from Proposition \ref{prop:1.1}. Then it follows from the definition of $\widehat\varphi$ in Proposition \ref{prop:2.1} that
\begin{equation*}
\widehat\varphi(db)=\sum_{(a)} \langle S(a_{(2)}),S\inv(d)\rangle \varepsilon(a_{(1)})=\langle a.d\rangle
\end{equation*}
so that indeed $a=\widehat\varphi(\,\cdot\,b)$.
\ebew

\iffalse %%%%%%%%% Old too long proof: 

\bew
Take $a'\in A$ and $b'\in B$. Then 
\begin{align*}
\langle a',b'b \rangle
&=\sum_{(a')} \langle a'_{(1)},b'\rangle\langle a'_{(2)},b\rangle \\
&=\sum_{(a')} \langle a'_{(1)},b'\rangle\psi(S(a'_{(2)})a) \\
&=\sum_{(a')} \langle S(a'_{(1)}),{\widehat S}^{-1}(b')\rangle\psi(S(a'_{(2)})a) \\
&=\sum_{(S(a'))} \langle (Sa')_{(2)}),{\widehat S}^{-1}(b')\rangle\psi((Sa')_{(1)})a) \\
&=\sum_{(a)} \langle S(a_{(2)}),{\widehat S}^{-1}(b')\rangle\psi((S(a')a_{(1)}) \\
&=\sum_{(a)} \langle a_{(2)}),b'\rangle\psi((S(a')a_{(1)}).
\end{align*}
We have used Equation (\ref{eqn:1.3}) of Propositon \ref{prop:1.1}
\begin{equation*}
S((\psi\ot\iota)(c\ot 1)\Delta(a))=(\psi\ot \iota)\Delta(c)(a\ot 1))
\end{equation*}
with $c=S(a')$.
\vskip 3pt
We see that 
$b'b=\sum_{(a)} \langle a_{(2)},b'\rangle\psi((S(\,\cdot\,)a_{(1)}).$
By the definition of $\widehat\varphi$ we find
\begin{equation*}
\widehat\varphi(b'b)=\sum_{(a)} \langle a_{(2)},b'\rangle \varepsilon(a_{(1)})=\langle a,b'\rangle
\end{equation*}
It follows that $\widehat\varphi(\,\cdot\,b)=a$ if $b=\psi(S(\,\cdot\,)a)$. This proves the result.
\ebew

\mycomment{\rood I have the feeling that this could be done simpler.
\ssnl}{}
\fi %%%%%%%%%%%%%%%%%%%%%

In fact, this is the same result as in the previous proposition, but with $A$ and $B$ interchanged.
\ssnl
\oldcomment{\rood 
Hier mag wel een extra referentie komen naar het origineel werk.
\ssnl
Iets nauwkeuriger formuleren ?
\ssnl
}

The more common choice for the Fourier transform is the map $a\mapsto \varphi(\,\cdot\,a)$ where $a\in A$ and where $\varphi$ is a left integral on $A$ (as in Proposition \ref{prop:1.4}).  We will work however with the second one for reasons that will be explained  later. The choice is related with the fact that we will use the right integral to construct the GNS-space in Section 2 of \cite{VD-part2}.%\ref{sect:posint},
%a point of discussion we mentioned already in the introduction. 
\ssnl
In Proposition 1.12 of \cite{VD-part2} we prove Plancherel's formula in the case of a positive integral on a multiplier Hopf $^*$-algebra. 
The results we prove here in Proposition \ref{prop:2.6} and Proposition \ref{prop:2.7} are in fact non-involutive versions of this formula.
\ssnl
For more details on the Fourier transform in the context of duality for algebraic quantum groups, we also refer to \cite{VD-ft}.
\ssnl
Further we will use $\mathcal F$ to denote this Fourier transform. So $\mathcal F(a)=\psi(S(\,\cdot\,)a)$ for $a\in A$. For this, we need to fix a right integral $\psi$ first. 
 
\oldcomment{\ssnl
Make a comment on Plancherel's formula, with reference to the second paper.
\ssnl}{}
\oldcomment{We also need a reference to \cite{VD-ft}}{}

 \section{\hspace{-17pt}. Duality for multiplier Hopf algebras} \label{sect:duality}% \input artikel3.tex %\newpage
 \oldcomment{\ssnl
 This part is corrected. 
 \ssnl}{}
In the two previous sections we have collected some properties of an algebraic quantum group and its dual. In this section, we proceed with the study of such a pair. Because some results are also true, and somewhat more natural, for more general pairs of multiplier Hopf algebras, we start with this case first.
\ssnl
The results  we formulate here are found in different earlier papers on the subject. 
The notion of a pairing of multiplier Hopf algebras has been considered in later papers as well. See e.g.\ Section 6 in \cite{Dr-VD-Z} and Section 4 in the survey paper \cite{VD-Zh2}. In fact, more recently, generalizations of such pairings have been studied for other related objects, see e.g. \cite{T-VD-W, La-VD4, La-VD5}.
\ssnl
Here we give a comprehensive, but also a more logical treatment of these results.
\ssnl
So, in the first place, we consider two regular multiplier Hopf  algebras $A$ and $B$ together with a non-degenerate pairing $\langle \,\cdot\,,\cdot\,\rangle:A\times B\to \mathbb C$.
\ssnl
As we have done before, we use the following convention, also in this section.
% \mycomment{It would make sense to formulate this remark earlier.} % We have done this in the introduction.

\opm % This is important and so we formulate it as a separate remark - Should we add something in the introduction?
We use $\Delta$, $\varepsilon$ and $S$ for the coproduct, the counit and the antipode on $A$ as well as on $B$. We mostly do the same for other objects associated with the multiplier Hopf algebras $A$ and $B$. On the other hand, we will systematically use symbols $a,a',\dots$ and $b,b'\dots$ for elements in $A$ and $B$ respectively. Then using the same symbols for the coproducts, the counits and the antipodes should not lead to any confusion.
\eopm
% We begin our discussion with the more general case. %We have mentioned this already.
\nl 
\bf Pairings of regular multiplier Hopf algebras \rm
\nl
We begin in this section by looking at  a more general pairing of two algebras. Recall the following definition (Definition 3.1  in \cite{T-VD-W}).

\defin\label{defin:1.3.2}
Consider a non-degenerate pairing of non-degenerate algebras $A$ and $B$. Assume that there exists {\it four actions}, a left and a right action of $A$ on $B$ and a left and a right action of $B$ on $A$, given by the formulas
\begin{align*}
\langle a', a \tr b \rangle &= \langle a'a ,b \rangle %\tussen
&\langle a' , b\tl a \rangle = \langle aa', b\rangle \\
\langle b\tr a, b'\rangle &= \langle a ,b'b \rangle %\tussen
&\langle a\tl b,b' \rangle = \langle a, bb' \rangle
\end{align*}
where $a,a'\in A$ and $b,b'\in B$.  If moreover these actions are unital we call the pairing an \emph{admissible pairing}.
\edefin

Recall that the action $(a,b)\mapsto a\tr b$ of $A$ on $B$ is called unital if $B$ is spanned by elements of the form $a\tr b$ with $a\in A$ and $b\in B$. Similarly for the other cases.

\opm\label{opm:3.3}
Actions of multiplier Hopf algebras where first studied in \cite{Dr-VD-Z}. 
In \cite{Dr-VD}, the existence of actions as in the previous definition is part of the axioms for a pairing of multiplier Hopf algebras. Also in later work, \cite{T-VD-W}, the existence of actions for a pairing of algebras is formulated as a condition. See e.g.\ Definition 3.1 in \cite{T-VD-W}. In the theory of quantum hypergroups under different appearances, the existence of such pairings for a pair of algebras plays an important role, see \cite{La-VD4, La-VD5}.
\eopm

Using the fact that these actions are all unital, it is possible to extend the pairing on $A\times B$ to $M(A)\times B$ and to $A\times M(B)$ provided the algebras have local units. This is found in Proposition 3.2 of \cite{T-VD-W}. We give the precise formulation below and we include a proof for completeness.

\prop\label{prop:1.3.4} 
Suppose that we have an admissible pairing $(a,b)\mapsto \langle a,b\rangle$ of algebras $A$ and $B$. If the algebras have local units we can extend the pairing to $M(A)\times B$ and to $A\times M(B)$ using the formulas  
\begin{align}
\langle m, a \tr b \rangle &= \langle ma ,b \rangle %\tussen
&\langle m , b\tl a \rangle & = \langle am, b\rangle \label{eqn:3.1}\\
\langle b\tr a, n\rangle &= \langle a ,nb \rangle% \tussen
&\langle a\tl b,n \rangle &= \langle a, bn \rangle\label{eqn:3.2}
\end{align}
where now $a\in A$, $m\in M(A)$, $b\in B$ and $n\in M(B)$. 
\eprop

\bew
i) First we define $\langle m,b\rangle$ for $m\in M(A)$ and $b\in B$. We use that any element $b\in B$ can be written as $\sum_i a_i\tr b_i$. On the other hand, suppose that we have elements $a_i$ and $b_i$ such that  $\sum_i a_i\tr b_i=0$. Then, with $e\in A$  satisfying $ema_i=ma_i$ for all $i$, we get
\begin{equation*}
\sum_i\langle ma_i,b_i\rangle=\sum_i\langle ema_i,b_i\rangle=\sum_i\langle em, a_i\tr b_i\rangle=0.
\end{equation*}
It follows from these two considerations that we can define a bilinear map on $M(A)\times B$ by $\langle m,a\tr b\rangle=\langle ma,b\rangle$.
\ssnl
ii) We claim that also $\langle m,b\tl a\rangle=\langle am,b\rangle$. To show this we take an element $e\in A$ satisfying $e\tr (b\tl a)=b\tl a$. Then we get
\begin{equation*}
\langle m,b\tl a\rangle=\langle me,b\tl a\rangle=\langle ame,b\rangle.
\end{equation*}
We can choose $e$ so that also $ame=am$ and we find that $\langle m,b\tl a\rangle=\langle am,b\rangle$.
\ssnl
iii) Finally, if already $m\in A$ we have $\langle m,a\tr b\rangle=\langle ma,b\rangle$ by the definition of the action in Definition  \ref{defin:1.3.2}. So we do get an extension of the pairing on $A\times B$ to $M(A)\times B$ satisfying the Equations \ref{eqn:3.1}.
\ssnl
iv) In a completely similar way, the the pairing is extended to $A\times M(B)$ satisfying Equations \ref{eqn:3.2}.
\ebew

 Remark  that it is in general not possible to extend the pairing to $M(A)\times M(B)$.
\ssnl
The above result also applies to  the induced pairings of tensor products.

\opm %{\blauw Move this ?}
This result appears in different forms in previous papers. First see a remark following the proof of Proposition 3.4 in 
\cite{De-VD2}. Also in Section 2 of \cite{VD-W}, results of this type are found.
\eopm

Now we can say what is a pairing of regular multiplier Hopf algebras.

\defin\label{defin:1.3.6}
Let $(A,\Delta)$ and $(B,\Delta)$ be regular multiplier Hopf algebras. Assume that we have an admissible pairing of $A$ with $B$. Then we call it a \emph{pairing of multiplier Hopf algebras} if
\begin{equation}
\langle \Delta(a), b\ot b' \rangle = \langle a ,bb' \rangle \tussenen
\langle a\ot a' , \Delta(b) \rangle = \langle aa', b \rangle \label{eqn:3.3a}
\end{equation}
for $a,a'\in A$ and $b,b'\in B$.
\edefin

We consider the tensor product pairing on $(A\ot A) \times (B\ot B)$ and its extensions to the multiplier algebras as obtained in Proposition \ref{prop:1.3.4}. This result is available in Section 2 of \cite{Dr-VD} in  a more complicated form.  But one easily sees that this is equivalent with the relation
\begin{equation*}
\langle \Delta(a)(1\ot a'),b\ot b'\rangle=\langle a\ot a',(b\ot 1)\Delta(b')\rangle
\end{equation*}
for all $a,a'\in A$ and $b,b'\in B$. In this form, for algebraic quantum hypergroups, it is found in \cite{De-VD2}.
\ssnl
Note the following important remark about conventions.

\opm\label{opm:3.5}
\oldcomment{As we mentioned already in the introduction,  I don't think we have mentioned this already - Check}{}
In the operator algebra approach to quantum groups, it is more common to flip the coproduct on the dual. This would mean e.g.\ that the first formula in Equation (\ref{eqn:3.1}) is still true but that the second would be replaced by\\ $\langle a\ot a' , \Delta(b) \rangle = \langle a'a, b \rangle $.
\vskip 3pt
In this paper however, we will systematically {\it stick to the algebraic convention} as in Definition \ref{defin:1.3.6}.
\eopm
 
For the counits on $A$ and on $B$ we find
\begin{equation} 
\varepsilon(a)=\langle a,1 \rangle
\tussenen \varepsilon(b)=\langle 1,b\rangle \label{eqn:3.2a}
\end{equation}
where $1$ denotes the identity in the multiplier algebra M(B) (in the first case) and in $M(A)$ (in the second case). For the antipodes we have the relation
\begin{equation}
\langle S(a),b \rangle = \langle a, S(b) \rangle\label{eqn:3.3}
\end{equation}
where again $a\in A$ and $b\in B$. Since we only consider regular multiplier Hopf algebras, the antipodes are bijections of $A$ and of $B$ respectively so that we do not have to use the extensions of the pairings for the formulas involving the antipodes.
\ssnl
These properties follow from the assumptions on the pairing as both the counit and the antipode are unique if they exist.

\nl
We can formulate more results on this level, but we get nicer formulations in the case of regular multiplier Hopf algebras with integrals. We consider this case in the next item.
\nl
\bf Dual pairs of regular multiplier Hopf algebras with integrals \rm
\nl
Now we will assume that the multiplier Hopf algebras $A$ and $B$ are \emph{algebraic quantum groups}. This means that they are \emph{regular} and that \emph {they have integrals}. 
\ssnl
We will use $\varphi$ for a left integral and $\psi$ for a right integral. Recall that integrals are unique up to a scalar. Here we also  use the same symbols for integrals on $A$ and on $B$. 
\ssnl
In what follows, we fix a left integral $\varphi$ on $A$.
\ssnl
\oldcomment{ Dit is een herhaling en is hier niet meer op zijn plaats.}

\defin
We let $B$ be the dual $\widehat A$ of $A$ and we define the pairing by  
\begin{equation*}
\langle a, b\rangle=\varphi(ac)
\end{equation*}
where $a\in A$  and $b=\varphi(\,\cdot\,c)$ for $c\in A$. 
\edefin

Recall that all elements of the form $\varphi(\,\cdot\,c)$, $\varphi(c\,\cdot\,)$, $\psi(\,\cdot\,c)$ and $\psi(c\,\cdot\,)$ belong to the dual $\widehat A$ and that any element in $\widehat A$ is of any of such forms. See a remark before Proposition 3.12 in \cite{VD-alg}
\ssnl
The following result can be found in \cite{Dr-VD}, see a remark following Definition 2.8 in that paper. However, the proof that is given there is rather short and there seems to be no place where a more detailed argument is found. For this reason we include a proof here.

\prop
The pairing of $A$ with $B$, as defined above, is a pairing of multiplier Hopf algebras.
% as defined in \cite{Dr-VD}.
\eprop
\bew
i) By definition any element $b$ of $B$ is of the form $\varphi(\,\cdot\,c)$ for some $c\in A$. Because $\varphi$ is faithful, the element $c$ is unique. This implies that $\langle a,b\rangle$ is well defined by the formula in the definition.
\ssnl
ii) If given $a\in A$ such that $\langle a,b\rangle=0$ for all $b$, then $\varphi(ac)=0$ for all $c$. This implies that $a=0$ because $\varphi$ is faithful. On the other hand, if given $b$ such that $\langle a,b\rangle=0$ for all $a$, it follows that $\varphi(ac)=0$ for all $a$. Again this implies that $c=0$ and hence $b=0$. Therefore the pairing is non-degenerate.
\ssnl
iii) Let $b=\varphi(\,\cdot\,c)$ and $a\in A$. We have, for all $a'\in A$ that
\begin{equation*}
\langle a'a,b\rangle=\varphi(a'ac)=\langle a',b'\rangle
\end{equation*}
where $b'=\varphi(\,\cdot\,ac)$. We see that $a\tr b$ exists in $B$ and that it is equal to $\varphi(\,\cdot\, ac)$. Because $A^2=A$ this action is unital. On the other hand, we have
\begin{equation*}
\langle aa',b\rangle=\varphi(aa'c)=\varphi(a'c\sigma(a))
\end{equation*}
and we see that $b\tl a$ exists in $B$ and that it is equal to $\varphi(\,\cdot\,c\sigma(a))$. Also this action is unital because $A^2=A$.
\ssnl
iv) For the actions of $B$ on $A$, we need the formula for the product in $B$. First we consider $b'b$ when $b=\varphi(c\,\cdot\,)$ Then we find
\begin{align*}
\langle a,b'b\rangle
&=\langle \Delta(a),b'\ot b\rangle\\
&=\langle (\iota\ot\varphi)((1\ot c)\Delta(a)),b'\rangle\\
&=\langle S( (\iota\ot\varphi)(\Delta(c)(1\ot a))),b'\rangle.
\end{align*}
%Using the Sweedler notation we find
%\begin{equation*}
%\langle a,b'b\rangle=\sum_{(c)}\langle S(c_{(1)},b'\rangle \varphi(c_{(2)}a)
%\end{equation*}
We see that $b\tr a$ exists in $A$ and that it is $S( (\iota\ot\varphi)(\Delta(c)(1\ot a)))$.
% $\sum_{(c)}S(c_{(1)})\varphi(c_{(2)}a)$. 
As $\Delta(A)(1\ot A)=A\ot A$ and because $S$ is bijective, we get that this action is unital.
\ssnl
On the other hand, we consider $bb'$ where now $b=\psi(\,\cdot\, c)$. Then we find
\begin{align*}
\langle a,bb'\rangle
&=\langle \Delta(a),b\ot b'\rangle\\
&=\langle(\psi\ot\iota)(\Delta(a)(c\ot 1)),b'\rangle\\
&=\langle S((\psi\ot\iota)(a\ot 1)\Delta(c)),b'\rangle.
\end{align*}
So $a\tl b$ exists in $A$ and is equal to
\begin{equation*}
S((\psi\ot\iota)(a\ot 1)\Delta(c)).
\end{equation*}
Because $(A\ot 1)\Delta(A)=A\ot A$ we get that also this action is unital.
\ebew

We have used the formulas from Proposition \ref{prop:1.1}.
\ssnl

\oldcomment{Perhaps some references for some formulas to the first sections?\ssnl}{}

Now we recall the following results. They have been obtained in Section 4 of \cite{De-VD1}, but in a slightly different setting. Therefore here we include arguments for completeness (and the convenience of the reader).
\ssnl
The following result  is found in Proposition 4.12 of \cite{De-VD1}. Here we give a direct proof.

\prop\label{prop:3.8} 
There is a unique element $V$ in $M(B\ot A)$ so that 
$$\langle V,a\ot b \rangle = \langle a,b \rangle$$
for all $a\in A$ and $b\in B$. We also have $(\iota \ot S)V$ and $(S\ot\iota)V$ are  in $M(B\ot A)$.
\eprop

\bew
We use  the {\it flipped} pairing on $B\times A$, the tensor product of the flipped pairing with the original pairing as well as the extension of the pairing to $M(B\ot A)\times (A\ot B)$. 
\ssnl
We claim that $V$ is given by the following formulas. We use the Sweedler notation.
\vskip 3pt
As a {\it left} multiplier it satisfies
\begin{equation} 
V(b\ot a)=\sum_{(c)}\varphi(c_{(2)}\,\cdot\,)\ot S(c_{(1)})a \label{eqn:3.4}
\end{equation}
where $a\in A$ and $b=\varphi(c\,\cdot\,)$ with $c\in A$. As a {\it right} multiplier we have
\begin{equation}
(b\ot a)V=\sum_{(c)}\psi(\,\cdot\,c_{(1)})\ot aS(c_{(2)})\label{eqn:3.5}
\end{equation}
for $a\in A$ and now $b=\psi(\,\cdot\,c)$ with $c\in A$.
\vskip 3pt
Indeed, in the first case we have, for $a'\in A$ and $b'\in B$, 
\begin{align*}
\langle\sum_{(c)}\varphi(c_{(2)}\,\cdot\,)\ot S(c_{(1)}), a'\ot b'\rangle
&=\sum_{(c)}\varphi(c_{(2)}a')\langle S(c_{(1)}), b'\rangle \\
&=\sum_{(a')}\varphi(ca'_{(2)})\langle a'_{(1)}, b'\rangle \\
&=\sum_{(a')}\langle a'_{(2)}\ot a'_{(1)}, b\ot b'\rangle \\
&=\langle a',b'b\rangle = \langle b\tr a',b'\rangle\\
&=\langle V, b\tr a',b'\rangle=\langle V(b\ot 1),a'\ot b'\rangle.
\end{align*}
We have used that 
$S((\iota\ot\varphi)(\Delta(c)(1\ot a')))=(\iota\ot\varphi)((1\ot c)\Delta(a'))$
for all $a',c\in A$ (Equation (\ref{eqn:1.2}) of Proposition \ref{prop:1.1}). 
\ssnl
This will provide the formula
\begin{equation*}
 V(b\ot 1)=\sum_{(c)}\varphi(c_{(2)}\,\cdot\,)\ot S(c_{(1)}).
\end{equation*}

This not only proves the formula (\ref{eqn:3.4}) but it also implies that $(1\ot a)V(b\ot 1)$ belongs to $B\ot A$.  %{\blauw Hier mag wel wat meer uitleg bij!}
\ssnl
Similarly we find for the second case
\begin{align*}
\langle \sum_{(c)}\psi(\,\cdot\,c_{(1)})\ot S(c_{(2)}),a'\ot b'\rangle
&=\sum_{(c)}\psi(a'c_{(1)}) \langle S(c_{(2)}),b'\rangle \\
&=\sum_{(a')}\psi(a'_{(1)}c) \langle a'_{(2)},b'\rangle \\
&=\sum_{(a')}\langle a'_{(1)}\ot  a'_{(2)},b\ot b'\rangle \\
&=\langle a',bb'\rangle=\langle a'\tl b,b'\rangle
\end{align*}
Here we used 
$S((\psi\ot\iota)((a'\ot 1)\Delta(c)))=(\psi\ot\iota)(\Delta(a')(c\ot 1))$
for all $a',c\in A$ (Equation (\ref{eqn:1.3}) of Proposition \ref{prop:1.1}).
This will give the formula
\begin{equation*}
 (b\ot 1)V=\sum_{(c)}\psi(\,\cdot\,c_{(1)})\ot S(c_{(2)})
\end{equation*}
and this proves the formula (\ref{eqn:3.5}).  Moreover we get that also $(b\ot 1)V(1\ot a)\in B\ot A$.
\ssnl
To complete the proof, we use that 
\begin{align*}
\langle (b_1\ot 1)(V(b_2\ot 1)),a'\ot b'\rangle&=\langle V(b_2\ot 1),(a'\tl b_1)\ot b'\rangle=\langle (b_2\tr a'\tl  b_1),b'\rangle\\
\langle ((b_1\ot 1)V)(b_2\ot 1),a'\ot b'\rangle&=\langle (b_1\ot 1)V,(b_2\tr a')\ot b'\rangle=\langle (b_2\tr a'\tl b_1), b'\rangle.
\end{align*}
So we have $(b_1\ot 1)(V(b_2\ot 1))=(b_1\ot 1)V)(b_2\ot 1)$ and the formulas (\ref{eqn:3.4}) and (\ref{eqn:3.5}) indeed define a multiplier $V$.
\oldcomment{Sufficient?}{}
\ebew

%{\blauw We have to rewrite these two arguments, and use the formulas we have in the previous result. Perhaps we have to include these formulas in the formulation of the proposition there so that we can refer to them here.
%\ssnl}

The result in Proposition \ref{prop:3.8} above is no longer true in general. Consider e.g.\ a pairing of infinite-dimensional Hopf algebras. In that case, the linear functional on $A\ot B$ defined by $a\ot b\mapsto \langle a,b\rangle$, in general,  can not be obtained by pairing of $a\ot b$ with an element of $B\ot A$. Remember that in this case, the algebras are unital so that the multiplier algebra $M(B\ot A)$ coincides with $B\ot A$. %{\blauw Ook hier wat meer uitleg?}
\nl
Now we translate the properties of the pairing in terms of the element $V$.
\ssnl
\prop
We have 
\begin{equation}(\iota\ot\Delta)V=V_{12}V_{13}
\tussenen
(\Delta\ot\iota)V=V_{13}V_{23}\label{eqn:3.6}
\end{equation}
where we are using the {\it leg numbering notation} as explained in the introduction. The first equation is valid in $M(B\ot A\ot A)$ and the second one in $M(B\ot B\ot A)$. We use the extension of the maps $\iota\ot\Delta$ and $\Delta\ot\iota$ to the multiplier algebra $M(B\ot A)$.
\eprop

These formulas are essentially obtained by rewriting the two formulas in Definition  \ref{defin:1.3.6}.
\ssnl
Reformulation of the formulas (\ref{eqn:3.2}) give
$(\varepsilon\ot\iota)V=1$ and $(\iota\ot\varepsilon)V=1$.
\ssnl
For the antipode we find the following. In the proof, we use $m$ for the multiplication map from $A\ot A$ to $A$.

\prop
The element $V$ is invertible in $M(B\ot A)$ and the inverse is given by
\begin{equation*}
 V^{-1}=(S\ot \iota)V=(\iota\ot S)V.
\end{equation*}
\eprop

\bew
We have seen in Proposition \ref{prop:3.8} that the elements $(S\ot \iota)V$ and $(\iota\ot S)V$ belong to $M(B\ot A)$.
\ssnl
If we apply $m(\iota\ot S)$ on the last two factors of the first equation in (\ref{eqn:3.6}), we find for the left hand side $(\iota\ot \varepsilon)V\ot 1$ and for the right hand side the product of $V$ and $(\iota\ot S)V$. As $(\iota\ot \varepsilon)V=1$ we find $V(\iota\ot S)V=1\ot 1$. On the other hand, if we apply $m(S\ot\iota)$ (again on the last two factors of the first equation in (\ref{eqn:3.6}), we see that the product of $(\iota\ot S)V$ and $V$ is also equal to $1\ot 1$. Therefore $V$ is invertible and its inverse is $(\iota\ot S)V$.
\vskip 3pt
We could do the same with the second formula in (\ref{eqn:3.6}) and obtain that the inverse is also $(S\ot\iota)V$. On the other hand, the equality $(S\ot \iota)V=(\iota\ot S)V$ is a reformulation of Equation (\ref{eqn:3.3}).
\ebew

In the case of a multiplier Hopf $^*$-algebra we obtain that $V$ is a unitary element. Indeed, a straightforward reformulation of the equations in \ref{eqn:1.4 } gives that 
$V^*=(S\ot\iota)V$.

\opm
Also for the more general pairings,  as in the beginning of this section, it is possible to consider similar formulas and to give them a meaning in a certain framework. See again Section 4 in \cite{De-VD1}. 
\vskip 3pt
The results in the case of a dual pair of algebraic quantum groups as treated above are much nicer. We will therefore stick to that case.
\eopm
\oldcomment{Moeten we dit hier niet aanvullen met een paar referenties? Ik denk dat dit nu ok is.}{}
In the next item, we study the algebra generated by $A$ and $B$ subject to the {\it Heisenberg commutation relations} as already formulated in the introduction.
\nl
\bf The Heisenberg algebra for a pair of multiplier Hopf algebras \rm
\nl
Let us again begin with the case of a non-degenerate pairing of two regular multiplier Hopf algebras $A$ and $B$ as in the beginning of this section. Most of the following results are found in Section 6 of \cite{Dr-VD-Z}.

\lem
Let $A$ act from the left on $A$ by multiplication. Consider also the left action $(b,a)\mapsto b\tr a$ of $B$ on $A$. Then we have the following commutation rules:
$$b\tr (ax)=\sum_{(a)(b)} \langle a_{(2)},b_{(1)}\rangle\, a_{(1)}(b_{(2)}\tr x)$$
where $a,x\in A$ and $b\in B$. 
\elem

On the right hand side, the element $b_{(2)}$ is covered. Indeed, as the action is unital, any element in $A$ is the linear span of elements of the form $b\tr a$ and the existence of local units in $B$ will give that, for all $x\in A$, there is an element $f\in B$ satisfying $f\tr x=x$. Successively, the element $a_{(1)}$ will be covered.
\oldcomment{Moeten we ook hier geen referentie hebben. Die hebben we nu al gegeven hierboven.}\zwart
\ssnl
This leads to the following proposition.

\prop
Consider the vector space $A\ot B$. It is an associative algebra for the product defined by
\begin{equation}
(a\ot b)(a'\ot b')=\sum_{(a')(b)} \langle a'_{(2)},b_{(1)}\rangle\,\, aa'_{(1)}\ot b_{(2)}b'.\label{eqn:3.7}
\end{equation}
It acts faithfully on $A$ from the left by 
$(a\ot b)\tr x=a(b\tr x)$.
\eprop

The proof is straightforward. It can be found in Section 6 of \cite{Dr-VD-Z}.
\ssnl
The commutation relations (\ref{eqn:3.7}) that determine the structure of the algebra are called the {\it Heisenberg commutation relations} and the algebra is the {\it Heisenberg algebra} associated with the pair $(A,B)$. 
\ssnl
We have the obvious non-degenerate embeddings $j_A$ and $j_B$ from $A$ and $B$ into $M(A\ot B)$ so that $a\ot b=j_A(a)j_B(b)$ when $a\in A$ and $b\in B$. If we identify $A$ and $B$ with their images in $M(A\ot B)$, the following notation makes sense. 

\notat 
In what follows we denote the Heisenberg algebra as $AB$ and we will also write the action of $AB$ on $A$ using the symbol $\tr$. We will systematically use letters $x,x',\dots$ for the elements in $A$ when they are acted upon by $AB$. % So we have $(ab)\tr x=\sum_..$
% NO this is not ok
\enotat 

Then the Heisenberg commutation rules are written as
\begin{equation*}
ba=\sum_{(a)(b)} \langle a_{(2)},b_{(1)}\rangle\,\, a_{(1)} b_{(2)}.
\end{equation*}
It is shown in \cite{Dr-VD} that the linear map $R$ from $A\ot B$ to $M(A\ot B)$, defined by 
$$R(a\ot b)= \sum_{(a)(b)} \langle a_{(2)},b_{(1)}\rangle\,\, a_{(1)}\ot b_{(2)},$$ 
is in fact a bijective map from $A\ot B$ to itself. This is  part of the set of axioms of a pairing of multiplier Hopf algebras. The inverse of this map is easily seen to satisfy 
\begin{equation}
 R\inv(a\ot b)=\sum_{(a)(b)} \langle S\inv(a_{(2)}),b_{(1)}\rangle\,\, a_{(1)}\ot b_{(2)}.\label{eqn:3.8}
\end{equation}

%%%%%%% Dit is nu niet meer nodig 

\iffalse %%%%%%%%%%%%%%%%%%
 This implies the following. It gives an equivalent form of the Heisenberg commutation rules.

\prop \label{prop:1.16}
In the Heisenberg algebra we have, for all $a\in A$ and $b\in B$,
\begin{equation*}
aS(b)=\sum_{(a),(b)}\langle a_{(2)},b_{(2)}\rangle S(b_{(1)})a_{(1)}.
\end{equation*}
\eprop

\bew
Using the formula (\ref{eqn:3.8}) for the inverse of $R$ we get
\begin{equation*}
ab= \sum_{(a)(b)} \langle S\inv(a_{(2)}),b_{(1)}\rangle\,\,b_{(2)}a_{(1)}
\end{equation*}
and if we replace $b$ by $S(b)$ we obtain
\begin{align*}
aS(b)&= \sum_{(a)(b)} \langle S\inv(a_{(2)}),S(b_{(2)})\rangle\,\,S(b_{(1)})a_{(1)}\\
&=\sum_{(a),(b)}\langle a_{(2)},b_{(2)}\rangle S(b_{(1)})a_{(1)}
\end{align*}
\ebew

\fi %%%%%%%%%%%%%

\iffalse %%%%%%%%
\bew
Let $a\in A$ and $b\in B$. In the Heisenberg algebra we find
\begin{align*}
\sum_{(a),(b)}\langle a_{(2)},b_{(2)}\rangle S(b_{(1)})a_{(1)}
&=\sum_{(a),(b),(S(b_{(1)}))}\langle a_{(3)}, b_{(2)} \rangle \langle a_{(2)}, S(b_{(1)})_{(1)} \rangle  a_{(1)} S(b_{(1)})_{(2)} \\
&=\sum_{(a),(b)}\langle a_{(3)}, b_{(3)} \rangle \langle a_{(2)}, S(b_{(2)}) \rangle  a_{(1)} S(b_{(1)}) \\
&=\sum_{(a),(b)}\langle a_{(2)},  S(b_{(2)})b_{(3)} \rangle  a_{(1)} S(b_{(1)}) \\
&=\sum_{(a)} \varepsilon(a_{(2)})   a_{(1)} S(b)=aS(b).
\end{align*}
\ebew
\fi %%%%%%%%%

%\begin{commentaar}
%Er is onduidelijkheid over het gebruik van $AB$. Is dat ok? Of gebruiken we toch eerder $C$?
%\ssnl
%Te beslissen! Waar komt dit weer vandaan?
%\ssnl
%\end{commentaar}

We now consider again the case of an {\it algebraic quantum group} $A$ and its {\it dual} $\widehat A$.
\nl
\bf The Heisenberg algebra for an algebraic quantum group \rm
\nl
The commutation rules determining the Heisenberg algebra now are equivalent (as we see from the proof below) with the basic formulas we prove in the following proposition.

\prop\label{prop:3.15}
For all $a\in A$ and $b\in B$ we have
$$\Delta(a)V=V(a\ot 1)
\tussenen
V\Delta(b)=(1\ot b)V.$$
The first formula holds in the multiplier algebra of $AB\ot A$ while the second one is true in the multiplier algebra of $B\ot AB$.
\eprop

\bew
To prove the first formula, we take the pairing with an element $b$ of $B$ in the second factor.
For the left hand side we get 
\begin{equation}
\sum_{(a)(b)}a_{(1)}b_{(2)}\langle a_{(2)},b_{(1)}\rangle\label{eqn:3.9}
\end{equation}
whereas for the right hand side we get $ba$. We use that $V$ is the duality in $M(B\ot A)$ so that pairing it with an element $b'$ of $B$ in the second factor precisely yields this element $b'$.
\ssnl
To prove the second formula, we take the pairing with an element $a$ of $A$ in the first factor. For the left hand side we find again (\ref{eqn:3.7}) and for the right hand side we get again $ba$.
\ebew

Also in the more general case, we can give a meaning to these formulas, very much along the line of the proofs. See e.g.\ \cite{VD-VK}.
\ssnl
\oldcomment{We have to add this here and comment on in terms of $V$ - Move this remark?}
In the case of an algebraic quantum group $A$  we have the following characterization of the Heisenberg algebra. See Proposition 6.7 in \cite{Dr-VD-Z}.

\prop\label{prop:3.16}%{prop:HsIso}
If $A$ is an algebraic quantum group, 
%and if $B$ is $\widehat A$, 
then the Heisenberg algebra is isomorphic with the algebra of operators on $A$ spanned by rank one operators of the form $x\mapsto a\langle x,b\rangle$ where $a\in A$ and $b\in B$.
\eprop

\bew
Let $a,x\in A$ and $b=\varphi(c\,\cdot\,)$ with $c\in A$. Then 
\begin{align*}
(ab)\tr x&=\sum_{(x)}ax_{(1)}\langle x_{(2)},b \rangle = \sum_{(x)}ax_{(1)}\varphi(cx_{(2)})\\
&=\sum_{(c)}aS(c_{(1)})\varphi(c_{(2)}x)=\sum_i a_i \langle x,b_i\rangle
\end{align*}
where we write $\sum_i a_i \ot b_i$ for $\sum_{(c)}aS(c_{(1)})\ot \varphi(c_{(2)}\,\cdot\,)$. On the other hand, we know that $(A\ot 1)(S\ot \iota)\Delta(A)=A\ot A$ for a regular multiplier Hopf algebras. Then the result follows. 
%{\blauw Reference?}
\ebew

We see that in this case, the Heisenberg algebra is only depending on the pairing of the two spaces $A$ and $B$. It is the algebra $A\diamond B$, with underlying vector space $A\ot B$ and product given by
 \begin{equation*}
 (a_1\ot b_1)\diamond (a_2\ot b_2)=\langle a_2,b_1\rangle a_1\ot b_2.
\end{equation*}
For this reason, we sometimes think of $AB$ as the Heisenberg algebra, {\it together} with the embeddings of $A$ and $B$ in the multiplier algebra $M(AB)$. The Heisenberg algebra alone forgets too much of the underlying structure.
\nl
\bf The action of $V$ on $A\ot A$ \rm
\nl
In what follows, we have the pairing of an algebraic quantum group with its dual. We consider again the left action of $B$ on $A$, as well as the left action of $A$ on itself, given by multiplication. This yields a left action of $B\ot A$ on $A\ot A$. This action is still unital and so it has an extension to a left action of the multiplier algebra $M(B\ot A)$ on $A\ot A$. 
\ssnl
We use this in the formulation of the following result. It is essentially obtained in Section 4 of \cite{De-VD1}.

\prop\label{prop:3.17}%{prop:Vact}
The element $V$ acts on $A\ot A$ as the canonical map $T:A\ot A\to A\ot A$, defined by $T(x\ot x')=\Delta(x)(1\ot x')$.
\eprop

\bew
We will use a Sweedler type notation $V=v_{(1)}\ot v_{(2)}$ in what follows.
For $x,x'\in A$ and $b\in B$ we find
\begin{align*}
V\tr(x\ot x')
&=(v_{(1)}\tr x) \ot v_{(2)}x'\\
&=\sum_{(x)}\langle x_{(2)} ,v_{(1)}\rangle x_{(1)}\ot v_{(2)}x' \\
&=\sum_{(x)} x_{(1)}\ot x_{(2)}x' = T(x\ot x') \\
\end{align*}

\vskip -15pt
\ebew
%\begin{align*}
%(v_{(1)}\tr x)\, \langle v_{(2)}x',b\rangle
%=\sum_{(x)}x_{(1)}\langle v_{(1)},x_{(2)}\rangle \langle v_{(2)}x',b\rangle\\
%&=\sum_{(x)}x_{(1)}\langle V(1\ot x'),x_{(2)}\ot b\rangle\\
%&=\sum_{(x)}x_{(1)}\langle x_{(2)},x'\tr b\rangle\\
%&=\sum_{(x)}x_{(1)}\langle x_{(2)}x',b\rangle.
%\end{align*}
%This proves that $(v_{(1)}\tr x)\ot v_{(2)}x'=\Delta(x)(1\ot x')$ for all $x,x'\in A$. 
%\ebew
%\begin{commentaar}
%Het lijkt erop  dat we hier in het bewijs een omweg maken?  Beter formuleren!
%\ssnl
%\end{commentaar}

Because 
$$T(ax\ot x')=\Delta(a)\Delta(x)(1\ot x')=\Delta(a)T(x\ot x')$$ 
for all $a,x,x'\in A$, we see again that $V(a\ot 1)=V\Delta(a)$, a formula that we proved already in the multiplier algebra $M(AB\ot A)$. Because $V$ is invertible we also have 
$\Delta(a)=V(a\ot 1)V^{-1}$ in $M(AB\ot A)$. If we combine this with the formulas in \ref{eqn:3.6}, we get the pentagon equation $V_{12}V_{13}V_{23}=V_{23}V_{12}$ for $V$. It holds in the multiplier algebra $M(B\ot AB\ot A)$.
\nl
\bf Application of the Fourier transform \rm
\nl
We will now finish this section by the formula we get for $V$ if we transform it with the Fourier transform. We have discussed the Fourier transform at the end of the previous section. Recall that we use $\mathcal F$ defined as $\mathcal F(x)=\psi(S(\,\cdot\,)x)$ for $x\in A$ where $\psi$ is a fixed right integral on $A$. 
%\ssnl
%It converts multiplication to convolution and vice versa. This is the content of the following proposition.

\ssnl
As expected, the Fourier transform converts multiplication operators to {\it convolution operators} and vice versa. This is the content of the following proposition.

\prop\label{prop:3.18}%{prop:conversion} 
For $a,x\in A$ and $b\in B$ we find
$$\mathcal F(ax)=\mathcal F(x)\tl S^{-1}(a)
\tussenen
\mathcal F(b\tr x)=b\mathcal F(x).$$
\eprop

\bew
i) Let $a,a',x\in A$. then
\begin{align*}
\langle a',\mathcal F(ax)\rangle
&=\psi(S(a')ax)=\psi(S(S^{-1}(a)a')x)\\
&=\langle S^{-1}(a)a',\mathcal F(x)\rangle
=\langle a',\mathcal F(x)\tl S^{-1}(a)\rangle
\end{align*}
ii) Let $a,x\in A$ and $b \in B$. Then
\begin{align*}
\langle a,\mathcal F(b\tr x)\rangle
=\psi(S(a)(b\tr x))
&=\sum_{(x)}\psi(S(a)x_{(1)})\langle b,x_{(2)}\rangle\\
&=\sum_{(S(a))}\psi(S(a)_{(1)}x)\langle b,S^{-1}(S(a)_{(2)})\rangle\\
&=\sum_{(a)}\psi(S(a_{(2)})x)\langle b,a_{(1)})\rangle\\
&=\langle a, b\mathcal F(x)\rangle
\end{align*}
\ebew

As a consequence of this result, we find how the duality $V$, acting on $A\ot A$ as the canonical map, shown in Proposition \ref{prop:3.17}, transforms. %We have included this at the end of Section 1, see Proposition \ref{prop:3.19}.
\ssnl
Using these properties we get the following result.

\prop\label{prop:3.19}
For $x,x'\in A$ we find
$$(\mathcal F\ot \mathcal F)(\Delta(x)(1\ot x'))=((S^{-1}\ot\iota)\Delta(y'))(y\ot 1)$$
where $y=\mathcal F(x)$ and $y'=\mathcal F(x')$.
\eprop

\bew
Let $x,x'\in A$ and write $y,y'$ for their Fourier transforms in $\widehat A$. Then using the Sweedler type notation for $V$ as we did before, we find
\begin{align*}
(\mathcal F\ot \mathcal F)V(x\ot x')
&= \mathcal F (v_{(1)}\tr x)\ot \mathcal F(v_{(2)}x') \\
&= v_{(1)}\mathcal F(x)\ot (\mathcal F(x')\tl S^{-1}(v_{(2)})) \\
&= v_{(1)}y\ot (y'\tl S^{-1}(v_{(2)})) \\
&=\sum_{(y')}v_{(1)}y\ot \langle y'_{(1)},S^{-1}(v_{(2)})\rangle y'_{(2)} \\
&=\sum_{(y')}v_{(1)}y\ot \langle S^{-1}(y'_{(1)}),v_{(2)}\rangle y'_{(2)} \\
&=\sum_{(y')}S^{-1}(y'_{(1)})y\ot y'_{(2)} \\
&=((S^{-1}\ot\iota)\Delta(y'))(y\ot 1).
\end{align*}
%{\blauw We use that we have the duality $V$}
%\ssnl
Remark that in all these expressions, we have the necessary coverings when using the Sweedler notation. In all but one of the cases, where we use the Sweedler type of notation for $V$, there is no problem as the two legs are covered. There is also no problem with the covering of the legs of $\Delta(y')$ in the last two expressions. Only when the Sweedler notation is used for both, one might need to multiply with an extra element in the second factor from the right.
\ebew 
\oldcomment{
Opmerking over het gebruik van de Sweedler type notatie voor $V$. Iets over de covering. Dit lijkt me nu ok.}{}
\ssnl

Remark that the map 
$$y\ot y'\mapsto ((S^{-1}\ot\iota)\Delta(y'))(y\ot 1)$$ 
we find above is the inverse of the map 
$y\ot y'\mapsto \Delta(y')(y\ot 1)$ from $\widehat A\ot \widehat A$ to itself.
%\nl
\oldcomment{Het volgende moet vervangen worden door een verwijzing naar part 2. We  kunnen dat ook schrijven in de conclusions section.\snl}{}
\opm
In the second paper on the subject \cite{VD-part2}, we work with a multiplier Hopf $^*$-algebra $A$ with positive integrals. This  allows to embed $A$ in a Hilbert space $\mathcal H$ and obtain a non-degenerate $^*$-representation of the Heisenberg algebra by means of bounded operators on $\mathcal H$. The duality $V$ now will act as a unitary operator. 
\ssnl
The Hilbert space is the GNS-space obtained from the right integral. The Fourier transform that we use above turns out to be a unitary from $\mathcal H$ to the GNS-space obtained from the left integral on the dual $\widehat A$. Then the map $y\ot y'\mapsto ((S^{-1}\ot\iota)\Delta(y'))(y\ot 1)$ is a unitary on this new Hilbert space, usually denoted as $W$,  and its inverse $W^*$ is the map $y\ot y'\mapsto \Delta(y')(y\ot 1)$. 
\eopm
\oldcomment{
Dit verklaart ook waarom we dit stukje hier invoeren! Dit mag misschien benadrukt worden.}{}

 \section{\hspace{-17pt}. Conclusions and further remarks} \label{sect:conclusions} %\input artikel4.tex %\newpage

In the \it first section \rm of this paper, we have collected some basic formulas relating the various objects that come with a multiplier Hopf algebra with integrals. In the second section, we have done the same for the objects of the dual as related with those of the original  multiplier Hopf algebra. We have considered the Fourier transforms. 
\ssnl
Finally, in the third section, we have studied the duality $V$, defined in the mulitplier algebra $M(\widehat A\ot A)$, where $A$ is a multiplier Hopf algebra with integrals and $\widehat A$ its dual. The duality $V$ acts on $A\ot A$ as the canonical map $a\ot a'\mapsto \Delta(a)(1\ot a')$. We also calculated the transformation of this map by means of the Fourier transform. It gives the inverse of the map $b\ot b'\mapsto \Delta(b')(b\ot 1)$.
\ssnl
These results are not new, but scattered in the literature.
\ssnl
In the second paper on this subject, \cite{VD-part2}, we treat the case of a multiplier Hopf $^*$-algebra with positive integrals. For this case, we consider an Hilbert space realization of the Heisenberg algebra by considering the GNS representation associated with a positive right integral. Then the duality $V$ becomes a unitary operator.
\ssnl
Among other things, we give a proof of the fact that the scaling constant $\tau$, defined by $\varphi\circ S^2=\tau\varphi$, is trivial and also that $\varphi\circ S$ is a positive right integral when $\varphi$ is a positive left integral.  
\ssnl
There are still a few open problems in the theory of multiplier Hopf algebras, in particular about the non-regular ones. For a regular multiplier Hopf algebra, the antipode $S$ maps $A$ to $A$ and is bijective. Because there are Hopf algebras with an antipode that is not bijective, there are multiplier Hopf algebras that are not regular. Then in principle, the antipode need not have range in $A$, but in the multiplier algebra $M(A)$. It should be possible to construct such examples.
\ssnl
Finally, also various types of pairings of multiplier Hopf algebras that are not Hopf algebras could be  interesting objects to study.

% \section{\hspace{-17pt}. } % \input artikel4.tex
% \section{\hspace{-17pt}. } % \input artikel5.tex

%%%%%%%%%% Appendices  %%%%%%%%%%%%%%%%%%%%%%%%%%

% Zaken aan te passen voor het geval we met appendices werken.

% Het volgende verandert de nummering van de sections naar hoofdletters voor de appendices

\renewcommand{\thesection}{\Alph{section}} 

\setcounter{section}{0}

% Het volgende dient voor een betere spatiëring als gevolg van een bredere A.1. bijvoorbeeld.

% \input heading-appendix.tex % Nieuwe instructies voor de appendices

\renewenvironment{stelling}{\begin{itemize}\item[ ]\hspace{-28pt}\bf Theorem \rm }{\end{itemize}}
\renewenvironment{propositie}{\begin{itemize}\item[ ]\hspace{-28pt}\bf Proposition \rm }{\end{itemize}}
\renewenvironment{lemma}{\begin{itemize}\item[ ]\hspace{-28pt}\bf Lemma \rm }{\end{itemize}}

%\section{\hspace{-17pt}. Appendix. ??}\label{sect:appA}  \input artikel8a.tex \newpage

%%%%%%%% References %%%%%%%%%%%%%

% \input artikel9.tex % Referenties

\end{document}